\documentclass[11pt]{article}
\usepackage{amscd,amssymb,epic,amsmath}

\usepackage[dvips]{graphicx}
\usepackage[mathscr]{eucal}
\usepackage[all]{xy}

\def\itemi{Section}
\def\itemii{Subsection}

\newtheorem{proposition}{Proposition}
\newtheorem{lemma}{Lemma}
\newtheorem{corollary}{Corollary}
\newtheorem{theorem}{Theorem}

\newtheorem{remark}{Remark}

\title{\Large \bf Plancherel measure for the quantum matrix ball - 1}

\author{O.Bershtein$^\dagger$, Ye.Kolisnyk}

\date{Institute for Low Temperature Physics and Engineering,\\ 47 Lenin ave.
61103, Kharkov, Ukraine.\\ e-mail: bershtein@ilt.kharkov.ua,
evgen.kolesnik@gmail.com}

\begin{document}

\maketitle

Abstract: The Plancherel formula is one of the celebrated results of
harmonic analysis on semisimple Lie groups and their homogeneous spaces. The
main goal of this work is to find a $q$-analogue of the Plancherel formula
for spherical transform on the unit matrix ball. Here we present an explicit
formula for the radial part of the Plancherel measure. The q-Jacobi
polynomials as spherical functions naturally arise on the way.

\makeatletter
\let\@thefnmark\relax
\@footnotetext{$^\dagger\ $ The first author was partially supported by the
N.~I.~Akhiezer fund.} \makeatother

\section{Introduction}

Let us recall one of the most common problems of harmonic analysis on
homogenous spaces. Let $G$ be a real Lie group, $K$ be a closed subgroup and
$d\nu$ be a $G$-invariant Haar measure on $X=K\backslash G$. The
representation of $G$ by right shifts in $L^2(X,d\nu)$
$$
R(g): f(x)\mapsto f(xg), \qquad x \in X, g \in G\
$$
is strongly continuous and unitary. It is called a quasiregular
representation. The problem is to find a decomposition of $R$ into
irreducible representations.

A special case of the Riemannian symmetric space $X=K\backslash G$ and its
isometry group was studied in detail (\cite{Helg}, p. 192;\ \cite{Helg1}, p.
506). Harmonic analysis for these spaces was developed by E.~Cartan,
I.~Gelfand, F.~Berezin, Harish-Chandra, S.~Gindikin, F.~Karpelevich.

The problem of harmonic analysis is closely connected with the following.
Consider the algebra $D_G(X)$ of all $G$-invariant differential operators on
$X$. An important result of the representation theory is that the
decomposition of $R$ can be obtained by using common eigenfunctions of
operators from $D_G(X)$. Namely, the shifts of a common eigenfunction
generate an irreducible subrepresentation of $R$.

In the case of a Riemannian symmetric space the algebra $D_G(X)$ is finitely
generated and commutative (\cite{Helg}, p. 431). Consider a set of the
generators of $D_G(X)$ and their restrictions $\mathcal{L}_1,\mathcal{L}_2,
\cdots , \mathcal{L}_r$ onto the subspace of smooth $K$-invariant functions
on $X$ with compact support. The Plancherel measure is a Borel measure in
$r$-dimensional space, and the problem is to find this measure.

Let us describe briefly how to find common eigenfunctions. Recall that an
irreducible strongly continuous unitary representation $T$ is called a
representation of type I if it contains a nonzero $K$-invariant vector $v$.
We can assume that $(v,v)=1$. The function $f(g)=(T(g)v,v)$ is called a
spherical function. It is constant on double cosets $K \backslash G/ K$, so
it corresponds to a $K$-invariant function on $X$. This function is a common
eigenfunction of operators $\mathcal{L}_1, \mathcal{L}_2, \cdots
,\mathcal{L}_r$. The problem of the decomposition of $K$-biinvariant
functions on $G$ in terms of the spherical functions naturally arises while
solving the general decomposition problem of $L^2(X, d\nu)$.

Consider a more special case. The homogeneous space $SU_{n,n} / S(U_n \times
U_n)$ is a Hermitian symmetric space of noncompact type. It has the standard
Harish-Chandra realization as the unit ball $$\mathbb{D} = \{ z \in
\mathrm{Mat}_n \; | \; ||z|| < I \}$$ (in the space of complex $n \times
n$-matrices with respect to the operator norm). It worth to be mentioned that
standard generators of $D_G(X)$ are well known and their common
eigenfunctions are Jacobi polynomials \cite{VD_CDO, Hoog_oper}.

Quantum bounded symmetric domains were introduced in 1998 by L. Vaksman and
S. Sinel'shchikov \cite{SV}. L. Vaksman with his collaborators managed to
develop the noncommutative complex analysis and representation theory on
quantum domains. A series of works were dedicated to quantum matrix balls
that are the simplest examples of quantum bounded symmetric domains
\cite{Vak01, Vaks, ShSinVak2001, Bershtein1, SSV4}.

In the case of quantum disk some problems of noncommutative harmonic
analysis are solved \cite{KorVak91, Kor}. In particular, explicit formulas
for the invariant integral, spherical functions, and the Plansherel measure
are obtained.


In this paper we generalize results mentioned above for the quantum matrix
ball case. Imitating the classical approach, we construct a family of
commuting 'q-differential' operators and find the exact formula for their
common eigenfunctions. We use spherical functions which appeare to be
{$q$-Jacobi} polynomials. We obtain the decomposition of the biinvariant
functions in terms of the spherical functions and the exact formula for the
'radial part' of the Plancherel measure.

The authors are grateful to L.~Vaksman for formulation the problems and
invaluable discussions. We also thank D.~Shklyarov for many important
remarks.

\section{The radial part the of $U_q \mathfrak{sl}_{2n}$-invariant integral}

\subsection{Preliminaries on the quantum matrix ball}

All results from the next two subsections form the basic notions in the theory of
quantum bounded symmetric domains. We refer to \cite{SV, SSV4} for the first
appearance and full consideration of these notions.

Let $q \in (0,1)$. All algebras are assumed to be associative and unital,
and $\mathbb C$ is the ground field.

Consider the well-known quantum universal enveloping algebra (see e.g.
\cite{KlSch}) $U_q\mathfrak{sl}_{2n}$ corresponding to the Lie algebra
$\mathfrak{sl}_{2n}$. Recall that $U_q\mathfrak{sl}_{2n}$ is a Hopf algebra
with the generators $\{E_i,\:F_i,\:K_i,\:K_i^{-1}\}_{i=1}^{2n-1}$ and the
relations

$$K_iK_j \,=\,K_jK_i,\quad K_iK_i^{-1}\,=\,K_i^{-1}K_i \,=\,1;$$
$$K_iE_i \,=\,q^{2}E_iK_i,\quad K_iF_i \,=\,q^{-2}F_iK_i;$$
$$K_iE_j \,=\,q^{-1}E_jK_i,\quad K_iF_j \,=\,qF_jK_i,\quad |i-j|=1;$$
$$K_iE_j \,=\,E_jK_i,\quad K_iF_j \,=\,F_jK_i,\quad |i-j| \,> \,1;$$
$$E_iF_j \,-\,F_jE_i \:=\:\delta_{ij}\frac{K_i-K_i^{-1}}{q-q^{-1}};$$
$$
E_i^2E_j\,-\,(q+q^{-1})E_iE_jE_i \,+\,E_jE_i^2 \:=\:0,\quad |i-j| \,=\,1;
$$
$$
F_i^2F_j \,-\,(q+q^{-1})F_iF_jF_i \,+\,F_jF_i^2\:=\:0, \quad |i-j|\,=\,1;
$$
$$E_iE_j-E_jE_i\,=\,F_iF_j-F_jF_i\,=\,0, \quad |i-j|\,> \,1.$$

The coproduct, the counit, and the antipode are defined as follows:
\begin{align*}
\triangle{E_j}&=E_j \otimes 1+K_j \otimes E_j,& \varepsilon(E_j)&=0,&
S(E_j)&=-K_j^{-1}E_j,\\ \triangle{F_j}&=F_j \otimes K_j^{-1}+1 \otimes F_j,&
\varepsilon(F_j)&=0,& S(F_j)&=-F_jK_j,\\ \triangle{K_j}&=K_j \otimes K_j,&
\varepsilon(K_j)&=1,& S(K_j)&=K_j^{-1},\qquad j=1,...,2n-1.
\end{align*}

Equip the Hopf algebra $U_q\mathfrak{sl}_{2n}$ with the involution $*$:
$$
(K_j^{\pm 1})^* = K_j^{\pm 1}, \quad E_j^* = \left\{
\begin{array}{rl}
   K_j F_j, & j \neq n, \\
  -K_j F_j, & j = n,
\end{array}
\right. \quad F_j^* = \left\{
\begin{array}{rl}
   E_j K_j^{-1}, & j \neq n, \\
  -E_j K_j^{-1}, & j = n.
\end{array}
\right.
$$

Then $U_q\mathfrak{su}_{n,n} \stackrel{\mathrm{def}}{=}
(U_q\mathfrak{sl}_{2n}, *)$ is an $*$-Hopf algebra. It is a quantum analogue
of the algebra ${U \mathfrak{su}_{n,n} \otimes_{\mathbb{R}} \mathbb{C}}$,
where $\mathfrak{su}_{n,n}$ stands for the Lie algebra of the noncompact
real Lie group $SU_{n,n}$.

Let $U_q \mathfrak{s} (\mathfrak{gl}_n \times \mathfrak{gl}_n) \subset U_q
\mathfrak{sl}_{2n}$ denote the Hopf subalgebra generated by $E_j, F_j, \, j
\neq n,$ and $K_i, K_i^{-1}, i=1,...,2n-1$. The corresponding Hopf
$*$-subalgebra in $U_q\mathfrak{su}_{n,n}$ is denoted by $U_q \mathfrak{s}
(\mathfrak{u}_n \times \mathfrak{u}_n)$.

Recall important definition of the weight module. A
$U_q\mathfrak{sl}_{2n}$-module $V$ is called a weight one if
\begin{equation*}
V=\bigoplus\limits_{\mathbf{\lambda}\in P}V_{\mathbf{\lambda}}, \qquad
V_\lambda=\left\{v\in V\;\left|\;K_iv=q^{\lambda_i}v,\quad
i=1,2,\ldots,2n-1\right.\right\},
\end{equation*}
where $\mathbf{\lambda}=(\lambda_1,\lambda_2,\ldots,\lambda_{2n-1})$ and $P$
is the weight lattice of the Lie algebra $\mathfrak{sl}_{2n}$. Nonzero
summand $V_{\mathbf{\lambda}}$ is called a weight subspace of the weight
$\lambda$.

Further, all $U_q \mathfrak{sl}_{2n}$-modules are assumed to be the weight
ones what allows us to introduce the linear operators $H_j$,
$j=1,\ldots,2n-1$, in $V$ such that
\begin{equation*}
H_jv=\eta_jv,\qquad v \in V_\eta.
\end{equation*}
Therefore, one can formally consider
\begin{equation*}
K_i^{\pm 1}=q^{\pm H_i}.
\end{equation*}

\bigskip

We recall a definition of the $*$-algebra $\mathrm{Pol}(\mathrm{Mat}_n)_q$
from \cite{SSV4}. First, let $\mathbb{C}[\mathrm{Mat}_n]_q$ denote the
well-known algebra with the generators $z_a^\alpha$, $a,\alpha=1,...,n,$ and
the relations
\begin{align}
& z_a^\alpha z_b^\beta-qz_b^\beta z_a^\alpha=0, & a=b \quad \& \quad
\alpha<\beta,& \quad \text{or}\quad a<b \quad \& \quad \alpha=\beta,
\label{zaa1}
\\ & z_a^\alpha z_b^\beta-z_b^\beta z_a^\alpha=0,& \alpha<\beta \quad
\&\quad a>b, & \label{zaa2}
\\ & z_a^\alpha z_b^\beta-z_b^\beta z_a^\alpha-(q-q^{-1})z_a^\beta
z_b^\alpha=0,& \alpha<\beta \quad \& \quad a<b.& \label{zaa3}
\end{align}
The algebra $\mathbb{C}[\mathrm{Mat}_n]_q$ is called the algebra of
holomorphic polynomials on the quantum $n \times n$ matrices space (see
\cite{KlSch}).

Similarly, let $\mathbb{C}[\overline{\mathrm{Mat}}_n]_q$ denote an algebra
with the generators $(z_a^\alpha)^*$, $a,\alpha=1,\dots,n$, and the
relations
\begin{flalign}
& (z_b^\beta)^*(z_a^\alpha)^* -q(z_a^\alpha)^*(z_b^\beta)^*=0, & a=b \quad
\& \quad \alpha<\beta, & \quad \text{or} \quad a<b \quad \& \quad
\alpha=\beta, \label{zaa1*}
\\ & (z_b^\beta)^*(z_a^\alpha)^*-(z_a^\alpha)^*(z_b^\beta)^*=0, &
\alpha<\beta \quad \& \quad a>b, & \label{zaa2*}
\\ & (z_b^\beta)^*(z_a^\alpha)^*-(z_a^\alpha)^*(z_b^\beta)^*-
(q-q^{-1})(z_b^\alpha)^*(z_a^\beta)^*=0,& \alpha<\beta \quad \& \quad a<b. &
\label{zaa3*}
\end{flalign}

Moreover, let $\mathbb C[\mathrm{Mat}_n \oplus \overline{\mathrm{Mat}}_n]_q$
denote an algebra with the generators $z_a^\alpha$, $(z_a^\alpha)^*$,
$a,\alpha=1,\dots,n$, relations (\ref{zaa1}) -- (\ref{zaa3*}), and
additional relations
\begin{equation*}
(z_b^\beta)^*z_a^\alpha = q^2 \sum\limits_{a',b'=1}^n \sum\limits_{\alpha',
\beta' = 1}^n R(b,a,b',a') R(\beta, \alpha, \beta', \alpha')
z_{a'}^{\alpha'} \left( z_{b'}^{\beta'} \right)^* + (1-q^2) \delta_{ab}
\delta^{\alpha \beta},
\end{equation*}
where $\delta_{ab}$, $\delta^{\alpha \beta}$ are Kronecker symbols and
$$
R(j,i,j',i') = \left\{
\begin{array}{cl}
  q^{-1}, &\quad i \neq j\ \&\ j=j'\ \&\ i=i', \\
  1, &\quad i=j=i'=j', \\
  -(q^{-2}-1), &\quad i=j\ \&\ i'=j'\ \&\ i'>i, \\
  0, &\quad \mbox{otherwise.}
\end{array}
\right.
$$

Finally, let $\mathrm{Pol}(\mathrm{Mat}_n)_q \stackrel{\rm def}{=} (\mathbb
C[\mathrm{Mat}_n \oplus \overline{\mathrm{Mat}}_n]_q,*)$ be an $*$-algebra
with the natural involution: ${*: z_a^\alpha \mapsto (z_a^\alpha)^*}$. The
algebra $\mathrm{Pol}(\mathrm{Mat}_n)_q$ is called the algebra of
polynomials on the quantum $n \times n$ matrices space (see \cite{KlSch}).

We now recall an irreducible $*$-representation of
$\mathrm{Pol}(\mathrm{Mat}_n)_q$ in a pre-Hilbert space. Let $\mathcal H$
denote the $\mathrm{Pol}(\mathrm{Mat}_n)_q$-module with one generator $v_0$
and the defining relations
$$(z_a^{\alpha})^*v_0=0, \quad a, \alpha=1,...,n.$$

Let $T_F$ denote the representation of $\mathrm{Pol}(\mathrm{Mat}_n)_q$
which corresponds to $\mathcal H.$ It is called the Fock representation. All
statements of the following proposition are proved in \cite{SSV4}.
\begin{proposition}\label{inner1}
\begin{enumerate}
\item $\mathcal {H}=\mathbb C[\mathrm{Mat}_n]_qv_0$.
 \item $\mathcal{H}$ is
a simple $\mathrm{Pol}(\mathrm{Mat}_n)_q$-module. \item There exists a
unique sesquilinear form $(\cdot,\cdot)$ on $\mathcal H$ with the following properties:\\
i) $(v_0,v_0)=1$; ii) $(fv,w)=(v,f^*w)$ for all $v,w \in \mathcal H$, $f \in
\mathrm{Pol}(\mathrm{Mat}_n)_q$. \item The form $(\cdot,\cdot)$ is positive
definite on $\mathcal H$.
\end{enumerate}
\end{proposition}

Also it is proved in \cite{SSV4} that $\mathrm{Pol}(\mathrm{Mat}_n)_q$ is a
$U_q \mathfrak{su}_{n,n}$-module algebra \footnote{i.e., the multiplication
in $\mathrm{Pol}(\mathrm{Mat}_n)_q$ is a morphism of $U_q
\mathfrak{su}_{n,n}$-modules, and the involutions in
$\mathrm{Pol}(\mathrm{Mat}_n)_q$ and $U_q \mathfrak{su}_{n,n}$ are
compatible.}. The action of generators of $U_q \mathfrak{su}_{n,n}$ is given
by the formulae
\begin{equation*}
H_nz_a^\alpha=\left \{\begin{array}{cl}2z_a^\alpha, &a=n \;\&\;\alpha=n,
\\ z_a^\alpha, &a=n \;\&\;\alpha \ne n \quad{\rm or}\quad a \ne n \;\&\;
\alpha=n, \\ 0, &{\rm otherwise,}\end{array}\right.
\end{equation*}
\begin{equation*}
F_nz_a^\alpha=q^{1/2}\cdot \left \{\begin{array}{cl}1, & a=n \;\&
\;\alpha=n, \\ 0, &{\rm otherwise,}\end{array}\right.
\end{equation*}
\begin{equation*}
E_nz_a^\alpha=-q^{1/2}\cdot \left \{\begin{array}{cl}q^{-1}z_a^nz_n^\alpha,
&a \ne n \;\&\;\alpha \ne n, \\
(z_n^n)^2, & a=n \;\&\;\alpha=n, \\ z_n^nz_a^{\alpha}, &{\rm
otherwise,}\end{array}\right.
\end{equation*}
for all ${a,\alpha=1,\ldots,n;}$ and with $k \ne n$
\begin{equation*}
H_kz_a^\alpha= \left \{\begin{array}{cl}z_a^\alpha ,& k<n \;\&\;a=k
\quad{\rm or}\quad k>n \;\&\;\alpha=2n-k, \\-z_a^\alpha ,& k<n \;\&\;a=k+1
\quad{\rm or}\quad k>n \;\&\;\alpha=2n-k+1, \\ 0 ,&{\rm
otherwise,}\end{array}\right.
\end{equation*}
\begin{equation*}
F_kz_a^\alpha=q^{1/2}\cdot \left \{\begin{array}{cl}z_{a+1}^\alpha ,& k<n
\;\&\;a=k, \\ z_a^{\alpha+1} ,& k>n \;\&\;\alpha=2n-k, \\ 0 ,&{\rm
otherwise,}\end{array}\right.
\end{equation*}
\begin{equation*}
E_kz_a^\alpha=q^{-1/2}\cdot \left \{\begin{array}{cl}z_{a-1}^\alpha ,& k<n
\;\&\; a=k+1, \\ z_a^{\alpha-1},& k>n \;\&\;\alpha=2n-k+1, \\ 0 ,& {\rm
otherwise.}\end{array}\right.
\end{equation*}

Let \begin{equation*}\Lambda_n = \{(\lambda_1, \lambda_2, \ldots, \lambda_n)
\in \mathbb{Z}_+^n \ | \ \lambda_1 \geq \lambda_2 \geq \ldots \geq
\lambda_n\}\end{equation*} be a set of partitions of length not larger than
$n$. Similarly to the classical case, one obtains the decomposition
$\mathbb{C} [\mathrm{Mat}_n]_q = \bigoplus_{\lambda \in \Lambda_n}
\mathbb{C} [\mathrm{Mat}_n]_{q, \bf{\lambda}}$ into a sum of $U_q
\mathfrak{s}(\mathfrak{u}_n \times \mathfrak{u}_n)$-isotypic components,
where $\mathbb{C} [\mathrm{Mat}_n]_{q, \bf{\lambda}}$ is a simple $U_q
\mathfrak{s}(\mathfrak{u}_n \times \mathfrak{u}_n)$-module with the highest
weight
\begin{equation*}
(\lambda_1-\lambda_2,...,\lambda_{n-1}-\lambda_n,2\lambda_n,
\lambda_{n-1}-\lambda_n,...,\lambda_1-\lambda_2).
\end{equation*}
This decomposition gives rise to the decomposition
\begin{equation*}
\mathcal{H}=\bigoplus_{\lambda \in \Lambda_n} \mathcal{H}_{\bf{\lambda}},
\qquad \mathcal{H}_{\bf{\lambda}}=\mathbb{C}
[\mathrm{Mat}_n]_{q,\bf{\lambda}}v_0.
\end{equation*}

Recall a definition of the quantum analogue of the Harish-Chandra embedding
of the Hermitian symmetric space $ S(U_n \times U_n) \backslash SU_{n,n}
\hookrightarrow \mathrm{Mat}_n$. Let $\mathbb{C}[SL_{2n}]_q$ denote the
well-known Hopf algebra with the generators $\{t_{ij}\}_{i,j=1,\ldots,2n}$
and the relations
\begin{flalign*}
& t_{\alpha a}t_{\beta b}-qt_{\beta b}t_{\alpha a}=0, & a=b \quad \& \quad
\alpha<\beta,& \quad \text{or}\quad a<b \quad \& \quad \alpha=\beta,
\\ & t_{\alpha a}t_{\beta b}-t_{\beta b}t_{\alpha a}=0,& \alpha<\beta \quad
\&\quad a>b,
\\ & t_{\alpha a}t_{\beta b}-t_{\beta b}t_{\alpha a}-(q-q^{-1})t_{\beta a}
t_{\alpha b}=0,& \alpha<\beta \quad \& \quad a<b,
\\ & \det \nolimits_q \mathbf{t}=1.
\end{flalign*}
Here $\det_q \mathbf{t}$ is the $q$-determinant of the matrix
$\mathbf{t}=(t_{ij})_{i,j=1,\ldots,2n}$ defined by
\begin{equation*}
\det \nolimits_q\mathbf{t}\stackrel{\mathrm{def}}{=}\sum_{s \in
S_{2n}}(-q)^{l(s)}t_{1\,s(1)}t_{2\,s(2)}\ldots t_{2n\,s(2n)},
\end{equation*}
with $l(s)=\mathrm{card}\{(i,j)|\;i<j \; \& \; s(i)>s(j) \}$. The
comultiplication $\Delta$, the counit $\varepsilon$, and the antipode $S$
are defined as follows:
$$
\Delta(t_{ij})=\sum_kt_{ik}\otimes t_{kj},\qquad
\varepsilon(t_{ij})=\delta_{ij},\qquad S(t_{ij})=(-q)^{i-j}\det \nolimits_q
\mathbf{t}_{ji},
$$
where $\mathbf{t}_{ji}$ is the matrix derived from $\mathbf{t}$ by
discarding its $j$-th row and its $i$-th column.

We equip $\mathbb C[SL_{2n}]_q$ with the standard $U_q
\mathfrak{sl}_{2n}$-module algebra structure as follows (see \cite{SSV4}): for
$k=1,..,2n-1$,
\begin{align}
E_k\cdot t_{ij}= q^{-1/2} & \begin{cases} t_{i\,j-1}, & k=j-1,
\\ 0, & \text{otherwise},
\end{cases} \qquad
F_k \cdot t_{ij}= q^{1/2}
\begin{cases}
t_{i\,j+1}, & k=j,
\\ 0, & \text{otherwise},
\end{cases}
\\ & K_k \cdot t_{ij}=
\begin{cases}
qt_{ij}, & k=j,
\\q^{-1}t_{ij}, & k=j-1,
\\ t_{ij}, & \text{otherwise}.
\end{cases}
\end{align}
Denote by $U_q \mathfrak{sl}^{\mathrm{op}}_{2n}$ the Hopf algebra obtained
from $U_q \mathfrak{sl}_{2n}$ by changing the multiplication to the opposite
one. We can also equip $\mathbb C[SL_{2n}]_q$ with a $U_q
\mathfrak{sl}^{\mathrm{op}}_{2n}$-module algebra structure as follows: for
$k=1,..,2n-1$,
$$
E_k \cdot t_{ij}= q^{-1/2}
\begin{cases}
t_{i+1\,j}, & k=i,
\\ 0, & \text{otherwise},
\end{cases} \qquad
F_k \cdot t_{ij}= q^{1/2}
\begin{cases}
t_{i-1\,j}, & k=i+1,
\\ 0, & \text{otherwise},
\end{cases}
$$
$$
K_k \cdot t_{ij}=
\begin{cases}
qt_{ij}, & k=i,
\\ q^{-1}t_{ij}, & k=i+1,
\\ t_{ij}, & \text{otherwise}.
\end{cases}
$$

So, $\mathbb C[SL_{2n}]_q$ is a $U_q \mathfrak{sl}^{\mathrm{op}}_{2n}
\otimes U_q \mathfrak{sl}_{2n}$-module algebra (see \cite{SSV4}). The
subalgebra
\begin{multline}\mathbb{C}[SL_{2n}]_q^{(U_q \mathfrak{s}(\mathfrak{gl}_n
\times \mathfrak{gl}_n))^{\mathrm{op}} \otimes U_q
\mathfrak{s}(\mathfrak{gl}_n \times \mathfrak{gl}_n)} = \{ f \in
\mathbb{C}[SL_{2n}]_q \quad | \\ \quad (\xi_1 \otimes \xi_2) f =
\varepsilon(\xi_1) \varepsilon(\xi_2) f, \qquad \xi_1 \in U_q \mathfrak{s}
(\mathfrak{gl}_n \times \mathfrak{gl}_n)^{\mathrm{op}}, \; \xi_2 \in U_q
\mathfrak{s} (\mathfrak{gl}_n \times \mathfrak{gl}_n)\}\end{multline} will
be referred as the subalgebra of $U_q \mathfrak{s}(\mathfrak{gl}_n \times
\mathfrak{gl}_n)$-biinvariants.

Equip $\mathbb{C}[SL_{2n}]_q$ with the involution given by
\begin{equation*}
t_{ij}^*= \mathrm{sign} [(i-n-1/2)(n-j+1/2)](-q)^{j-i}\det \nolimits_q
\mathbf{t}_{ij}.
\end{equation*}
It can be proved that $\mathbb{C}[w_0 SU_{n,n}]_q\stackrel{\mathrm{def}} {=}
(\mathbb{C}[SL_{2n}]_q,*)$ is an $U_q \mathfrak{su}_{n,n}$-module
$*$-algebra. It is a $q$-analogue of the algebra of regular functions on the
real affine algebraic manifold $w_0 SU_{n,n}$, where\footnote{The matrix
$w_0$ corresponds to the longest element of the Weyl group of the Lie
algebra $\mathfrak{sl}_{2n}$.} $$w_0 = \left(%
\begin{array}{cc}
  0 & -J \\
  J & 0 \\
\end{array}%
\right),  \qquad  J = \left(%
\begin{array}{ccccc}
  0 & 0 & ... & 0 & 1 \\
    0 & 0 &...& 1 & 0 \\
  & & ... & & \\
  0 & 1 & ... & 0 & 0 \\
    1 & 0 & ...& 0 & 0 \\
\end{array}%
\right).$$

For any multiindices $I=\{1 \le i_1<i_2<\ldots<i_k \le 2n\}$ and $J=\{1 \le
j_1<j_2<\ldots<j_k\le 2n\}$ we use the following standard notation for the
corresponding $q$-minor of the matrix $\mathbf{t}$:
$$
t_{IJ}^{\wedge k}\stackrel{\mathrm{def}}{=}\sum_{s \in
S_k}(-q)^{l(s)}t_{i_1j_{s(1)}} t_{i_2j_{s(2)}} \ldots  t_{i_kj_{s(k)}}.
$$
We now introduce short notation for the elements
\begin{equation}\label{t_x}
t=t_{\{1,2,\dots,n \}\{n+1,n+2,\dots,2n \}}^{\wedge n},\qquad x=tt^*.
\end{equation}
Note that $t,$ $t^*,$ and $x$ quasi-commute with all generators $t_{ij}$ of
$\mathbb{C}[SL_{2n}]_q$, and that $\mathbb{C}[w_0 SU_{n,n}]_q$ is an
integral domain (see \cite{Jo}). Let $\mathbb{C}[w_0 SU_{n,n}]_{q,x}$ be the
localization of $\mathbb{C}[w_0 SU_{n,n}]_q$ with respect to the
multiplicative set $x^{\mathbb Z_+}$ (see \cite{BrGood}). The following
statements are proved in \cite{SSV4}.

\begin{proposition}
There exists a unique extension of the $U_q \mathfrak{su}_{n,n}$-module
$*$-algebra structure from $\mathbb{C}[w_0 SU_{n,n}]_q$ to $\mathbb{C}[w_0
SU_{n,n}]_{q,x}$.
\end{proposition}

\begin{proposition}
There exists a unique embedding of the $U_q \mathfrak{su}_{n,n}$-module
$*$-algebras $$i:\mathrm{Pol}(\mathrm{Mat}_n)_q \hookrightarrow
\mathbb{C}[w_0 SU_{n,n}]_{q,x}$$ such that
\begin{equation*}
i(z_a^\alpha)= t^{-1}t_{\{1,2,\dots,n \}J_{a \alpha}}^{\wedge n},
\end{equation*}
where $J_{a \alpha}=\{n+1,n+2,\dots,2n \}\setminus \{2n+1-\alpha \}\cup \{a
\}$.
\end{proposition}

The last proposition gives us a $q$-analogue of the Harish-Chandra
embedding. It allows us to identify $\mathrm{Pol}(\mathrm{Mat}_n)_q$ with
its image in $\mathbb{C}[w_0 SU_{n,n}]_{q,x}$.

\subsection{The algebra of finite functions and the invariant
integral}\label{finite_fun_ideal}

 It is well known that in the classical case $q=1$  a positive definite ${\rm
SU}_{n,n}$-invariant integral can not be defined on the polynomial algebra
in the unit ball $\mathbb{D} \hookrightarrow \mathrm{Mat}_n$. However, it is
well defined on the space of smooth functions with compact support on
$\mathbb{D}$. These observations are still applicable for the quantum case.
Here we provide the definition and some basic properties of a $q$-analogue
of the algebra of finite functions following \cite{SSV_CJP}.

Let us consider a $U_q\mathfrak{su}_{n,n}$-module $*$-algebra ${\rm
Fun}(\mathbb{D})_q$ obtained from ${\rm Pol}({\rm Mat}_n)_q$ by adding a
genera\-tor $f_0$ and the relations
$$f_0=f_0^2=f_0^*,$$ $$(z_a^\alpha)^* f_0=0, \quad f_0 z_a^\alpha=0,
\qquad a,\alpha = 1,2,\ldots,n.$$ The $U_q \mathfrak{su}_{n,n}$-module
algebra structure can be extended from $\mathrm{Pol}(\mathrm{Mat}_n)_q$ to
${\rm Fun}(\mathbb{D})_q$ as follows:
\begin{equation*}\label{sl_f0_1}
H_n f_0 = 0, \qquad F_n f_0 = -\frac{q^{1/2}}{q^{-2}-1} f_0 (z^n_n)^*,
\qquad E_n f_0 = -\frac{q^{1/2}}{1-q^2} z^n_n f_0,
\end{equation*}
\begin{equation*}\label{sl_f0_2}
H_k f_0 = F_k f_0 = E_k f_0 = 0, \qquad k \neq n.
\end{equation*}
The two-sided ideal $\mathscr{D}(\mathbb{D})_q = {\rm Pol}({\rm Mat}_n)_q
f_0 {\rm Pol}({\rm Mat}_n)_q$ is a $U_q \mathfrak{su}_{n,n}$-module
$*$-subalgebra (see \cite{SSV_CJP}). The elements of the two-sided ideal
$\mathscr{D}(\mathbb{D})_q$ will be called finite functions on the quantum
matrix ball $\mathbb{D}$.

The Fock representation $T_F$ of ${\rm Pol}({\rm Mat}_n)_q$ can be extended
up to the representation of ${\rm Fun}(\mathbb{D})_q$, and so for every finite
function $f\in\mathscr{D}(\mathbb{D})_q$ there exists an operator $T_F(f)$, and
$$ T_F (\mathscr{D}(\mathbb{D})_q) = \{ A \in \mathrm{End} (\mathcal{H}) \;
| \; A|_{\mathcal{H}_{\lambda}}\neq 0 \; \text{for a finite set of indices}
\; \lambda \in \Lambda_n\}.$$

Consider the gradings
$$\mathbb{C}[\mathrm{Mat}_n]_{q,k} = \bigoplus_{|\lambda| = k}
\mathbb{C}[\mathrm{Mat}_n]_{q,\lambda}, \qquad k \in \mathbb{Z}_+,$$ and
$$\mathbb{C}[\overline{\mathrm{Mat}_n}]_{q,-k} = \bigoplus_{|\lambda| = k}
\mathbb{C}[\overline{\mathrm{Mat}_n}]_{q,\lambda}, \qquad k \in
\mathbb{Z}_+,$$ where $|\lambda| = \lambda_1 + \lambda_2 + \ldots +
\lambda_n$.

It is evident that
\begin{lemma}\label{extension_Fock}
The Fock representation $T_F$ has a unique extension to a representation of
the $*$-algebra $\mathrm{Fun}(\mathbb{D})_q$ such that the element $f_0$
maps to the orthogonal projection onto the vacuum subspace.
\end{lemma}

Let us keep the same notation $T_F$ for this extension.

\begin{proposition}\label{bijection}
The representation $T_F$ provides the isomorphism of the $*$-algebra
$\mathscr{D}(\mathbb{D})_q$ and the $*$-algebra of all finite\footnote{A linear
operator $A$ in $\mathcal{H}$ is called finite if $A \mathcal{H}_j = 0$ for all $j
\in \mathbb{Z}_+$ except a finite set.} linear operators in $\mathcal{H}$.
\end{proposition}

{\bf Proof.} $T_F$ is an $*$-representation. So, we have to prove that the
restriction of $T_F$ on $\mathscr{D}(\mathbb{D})_q$ is a bijective mapping
from $\mathscr{D}(\mathbb{D})_q$ to the algebra of all finite linear
operators in $\mathcal{H}$.

Let $\mathscr{D}(\mathbb{D})_{q,i,j} = \mathbb{C}[\mathrm{Mat}_n]_{q,i} \cdot
f_0 \cdot \mathbb{C}[\overline{\mathrm{Mat}_n}]_{q,-j}$. If $f \in
\mathscr{D}(\mathbb{D})_{q,i,j}$ then the linear operator $T_F(f)$, $f \in
\mathscr{D}(\mathbb{D})_{q,i,j}$ maps $\mathcal{H}_j$ to $\mathcal{H}_i$ and it
is equal to zero on $\bigoplus\limits_{k \neq j} \mathcal{H}_k$. We obtain a
linear mapping from $\mathscr{D}(\mathbb{D})_{q,i,j}$ to
$\mathrm{Hom}(\mathcal{H}_j, \mathcal{H}_i)$. It is surjective by Proposition
\ref{inner1}, and $$ \dim \mathscr{D}(\mathbb{D})_{q,i,j} =
\dim\,\mathrm{Hom}(\mathcal{H}_j, \mathcal{H}_i).$$ Thus, the representation
$T_F$ provides the isomorphism
\begin{equation*}\label{iso_finite}
  \mathscr{D}(\mathbb{D})_{q,i,j} = \mathbb{C}[\mathrm{Mat}_n]_{q,i}\; f_0\;
\mathbb{C}[\overline{\mathrm{Mat}_n}]_{q,-j} \cong \mathrm{Hom}(\mathcal{H}_j,
\mathcal{H}_i).
\end{equation*}
But $\mathscr{D}(\mathbb{D})_q = \bigoplus\limits_{i,j=0}^\infty
\mathscr{D}(\mathbb{D})_{q,i,j}$, and $\bigoplus\limits_{i,j=0}^\infty
   \mathrm{Hom}(\mathcal{H}_j, \mathcal{H}_i)$ in
$\mathrm{End}\,\mathcal{H}$ is the vector space of finite linear operators.
\hfill $\square$


\begin{proposition}\label{extandd2}
The representation $T_F$ provides the bijection of the space of $U_q
\mathfrak{s}(\mathfrak{gl}_n \times \mathfrak{gl}_n)$-invariants in
$\mathscr{D}(\mathbb{D})_q$ and the space of finite linear operators in
$\mathcal{H}$ that are scalars on every $U_q \mathfrak{s}(\mathfrak{gl}_n
\times \mathfrak{gl}_n)$-isotypic component $\mathcal{H}_\lambda$,
$\lambda \in \Lambda_n$.
\end{proposition}
 {\bf Proof.}

i) If $f$ is a $U_q \mathfrak{s}(\mathfrak{gl}_n \times
\mathfrak{gl}_n)$-invariant vector, then $T_F(f)$ maps a highest vector of
$\mathcal{H}_\lambda$ to a highest vector of a $U_q
\mathfrak{s}(\mathfrak{gl}_n \times \mathfrak{gl}_n)$-isotypic component with
the same weight.

ii) The action of $U_q \mathfrak{s}(\mathfrak{gl}_n \times \mathfrak{gl}_n)$
in $\mathcal{H}$ is multiplicity free.

iii) Now i) and ii) imply that if $f$ is a $U_q \mathfrak{s}(\mathfrak{gl}_n \times
\mathfrak{gl}_n)$-invariant vector, then $T_F(f)|_{\mathcal{H}_\lambda}$ is an
endomorphism of the simple $U_q \mathfrak{s}(\mathfrak{gl}_n \times
\mathfrak{gl}_n)$-module $\mathcal{H}_\lambda$. So $T_F(f)$ is scalar on
$\mathcal{H}_\lambda$, $\lambda \in \Lambda_n$. \hfill $\Box$

Denote the space of $U_q \mathfrak{s}(\mathfrak{gl}_n \times
\mathfrak{gl}_n)$-invariants in $\mathscr{D}(\mathbb{D})_q$ by
$$(\mathscr{D}(\mathbb{D})_q)^{U_q \mathfrak{s} (\mathfrak{gl}_n \times
\mathfrak{gl}_n)} = \{ f \in \mathscr{D}(\mathbb{D})_q \quad | \quad \xi f =
\varepsilon(\xi) f,\quad \xi \in U_q \mathfrak{s} (\mathfrak{gl}_n \times
\mathfrak{gl}_n)\}.$$

Denote $$\check{\rho}= \frac{1}{2} \sum\limits_{j=1}^{2n-1} j(2n-j) H_j.$$

The following proposition is also stated in \cite{SSV_CJP}.
\begin{proposition}\label{positiveIntegral}
The linear functional
\begin{equation}\label{invin_d}
\int f d\nu = (1-q^2)^{n^2} \mathrm{tr}(T_F(f)q^{-2\check{\rho}}),\qquad f
\in \mathscr{D}(\mathbb{D})_q,
\end{equation}
is a positive definite $U_q \mathfrak{sl}_{2n}$-invariant integral on
$\mathscr{D}(\mathbb{D})_q$, i.e.
$$ \int \xi f d\nu = \varepsilon(\xi) \int f d\nu, \qquad \xi \in U_q \mathfrak{sl}_{2n}$$
and
$$ \int f^* f d\nu > 0, \qquad \mathrm{for} \; f \neq 0.$$
\end{proposition}
For the sketch of the proof refer to \cite[\S 5]{SSV_FT}.

Further we consider the restriction of invariant integral \eqref{invin_d} to the
space of $U_q \mathfrak{s} (\mathfrak{gl}_n \times \mathfrak{gl}_n)$-invariants
in $\mathscr{D}(\mathbb{D})_q$. We will call this restriction the radial part.

\subsection{The radial part of the invariant integral}\label{x_spectr_integral}

In this \itemii $\;$ we will describe the support of radial part of the
invariant measure $d\nu$ and find an exact formula for the radial part of
the invariant integral.

Consider the elements of $\mathbb{C}[w_0 SU_{n,n}]_q$:
$$x_k\,=\,q^{k(k-1)} \sum_{\substack{I \subset \{1,2,\ldots,n\},\, J\subset
\{n+1,n+2,\ldots,2n\} \\ \mathrm{card} (I)= \mathrm{card} (J)=k}} q^{-2
\sum\limits_{m=1}^k (n-i_m)} (-q)^{\sum\limits_{m=1}^k (j_m-i_m-n)}\,
t_{I\,J}^{\wedge k} \, t_{I^c\,J^c}^{\wedge (2n-k)}.$$ It follows from the results of
\cite{BKV} that $x_k$, $k=1,2,\ldots,n$, are pairwise commuting self-adjoint
${U_q \mathfrak{s} (\mathfrak{u}_n \times \mathfrak{u}_n)}$-biinvariants. These
elements generate the subalgebra of all $U_q \mathfrak{s} (\mathfrak{u}_n
\times \mathfrak{u}_n)$-biinvariant elements in $\mathbb{C}[w_0 SU_{n,n}]$, as
follows from the results of \cite{BerS} and \cite{BKV}. So, $$\mathbb{C}[w_0
SU_{n,n}]^{(U_q \mathfrak{s}(\mathfrak{u}_n \times
\mathfrak{u}_n))^{\mathrm{op}} \otimes U_q \mathfrak{s}(\mathfrak{u}_n \times
\mathfrak{u}_n)}_q \cong \mathbb{C}[x_1, x_2, \ldots, x_n].$$

Denote by $T$ the $*$-representation of the $*$-algebra $\mathbb{C}[w_0
SU_{n,n}]_q$ corresponding to the permutation
$$\left(\begin{tabular}{cccccccccc}
  1  & 2  & \ldots &  n-1 & n  & n+1 & n+2  & \ldots &  2n-1   & 2n \\
  n+1 & n+2 & \ldots & 2n-1 & 2n & 1 & 2 & \ldots & n-1 & n \\
\end{tabular}\right)$$ (see \cite[\S
4]{SSV4}). This representation allows a unique extension to the
representation of $\mathbb{C}[w_0 SU_{n,n}]_{q,x}$, where $x$ is defined in
\eqref{t_x}. It is proved in \cite[\S 4,5]{SSV4} that the representation
$T_F$ is unitary equivalent to the restriction of the representation $T$ to
$\mathrm{Pol}(\mathrm{Mat}_n)_q$. Consider the short notation
$$
q^{\mu} = ( q^{\mu_1}, q^{\mu_2}, \ldots, q^{\mu_n}) \in \mathbb{C}^n,
\qquad \mu \in \mathbb{C}^n.
$$
It is also proved in \cite{BKV} that
$$ T(x_k)|_{\mathcal{H}_\lambda}  =
q^{k(k-1)} e_k(q^{-2(\lambda +\delta)}) ,\qquad k = 1,2,\ldots,n,\quad
\lambda \in \Lambda_n,$$ where $\delta = (n-1, n-2, \ldots, 1, 0) \in
\Lambda_n$, $e_k$ is the elementary symmetric polynomial in $n$ variables of
the degree $k$. So, the set of common eigenvalues of the operators $T(x_1)$,
$T(x_2)$, $\ldots$, $T(x_n)$ is
\begin{equation*}
\Sigma_{\mathbb{D}} = \{ ( e_1(q^{-2(\lambda+\delta)}), q^2
e_2(q^{-2(\lambda+\delta)}), \ldots, q^{n(n-1)} e_n(q^{-2(\lambda+\delta)})
) \; | \; \lambda \in \Lambda_n\}.
\end{equation*}

Thus the algebra $\mathbb{C}[w_0 SU_{n,n}]^{(U_q \mathfrak{s}(\mathfrak{u}_n
\times \mathfrak{u}_n))^{\mathrm{op}} \otimes U_q \mathfrak{s}(\mathfrak{u}_n
\times \mathfrak{u}_n)}_q$ can be identified with the algebra of polynomial
functions on $\Sigma_{\mathbb{D}}$. Following Propositions \ref{bijection},
\ref{extandd2} the algebra $\mathscr{D}(\mathbb{D})^{U_q \mathfrak{s}
(\mathfrak{gl}_n \times \mathfrak{gl}_n)}_q$ can be identified with the algebra
$\mathscr{D}(\Sigma_{\mathbb{D}})$ of functions $f(x_1,x_2,\ldots,x_n)$ with
finite support on $\Sigma_{\mathbb{D}}$.

\begin{lemma}\label{spectr_bij}
The mapping $$\Lambda_n \to \Sigma_{\mathbb{D}},$$
$$\lambda = (\lambda_1,\lambda_2,\ldots,\lambda_n) \mapsto ( e_1(q^{-2(\lambda+\delta)}), q^2
e_2(q^{-2(\lambda+\delta)}), \ldots, q^{n(n-1)} e_n(q^{-2(\lambda+\delta)})
)$$ is a bijection.
\end{lemma}
{\bf Proof.} The surjectivity follows from the definition of
 $\Sigma_{\mathbb{D}}$. Let us prove the injectivity. The function $q^{-l}$ is strictly
increasing as $ l \in [0, +\infty)$, so the mapping
$$\Lambda_n \to \mathbb{R}^n, \qquad \lambda \mapsto
q^{-2(\lambda+\delta)}$$ is an injection. Due to the Viet theorem, the
mapping
$$q^{-2(\lambda+\delta)} \mapsto ( e_1(q^{-2(\lambda+\delta)}), q^2
e_2(q^{-2(\lambda+\delta)}), \ldots, q^{n(n-1)} e_n(q^{-2(\lambda+\delta)})
)$$ is also an injection since $q^{-2(\lambda_1+n-1)} >
q^{-2(\lambda_2+n-2)} > \ldots > q^{-2\lambda_n}$ for any $\lambda \in
\Lambda_n$. Now we have the injectivity of the composition.
 \hfill $\Box$

Consider the algebra $\mathbb{C}[u_1,u_2,\ldots,u_n]$ and the injection
\begin{equation}\label{x_e_link_inj}
\mathbb{C}[x_1, x_2, \ldots, x_n] \hookrightarrow
\mathbb{C}[u_1,u_2,\ldots,u_n], \qquad
 x_k \mapsto q^{k(k-1)} e_k(u_1, \ldots, u_n).
\end{equation}
This injection allows one to identify the subalgebra $\mathbb{C}[w_0
SU_{n,n}]^{(U_q \mathfrak{s}(\mathfrak{u}_n \times
\mathfrak{u}_n))^{\mathrm{op}} \otimes U_q \mathfrak{s}(\mathfrak{u}_n \times
\mathfrak{u}_n)}_q$ with the algebra of all symmetric polynomials in variables
$u_1,u_2,\ldots,u_n$.

Specify $$\Delta_{\mathbb{D}} = \{ (q^{-2(\lambda+\delta)}\;|\; \lambda \in
\Lambda_n \}.$$ Let also $\mathscr{D}(\Delta_{\mathbb{D}})$ be the algebra of
functions $f(u_1,u_2,\ldots,u_n)$ with finite support on the set
$\Delta_{\mathbb{D}}$. Then
$$\mathscr{D}(\Delta_{\mathbb{D}}) \cong
 \mathscr{D}(\Sigma_{\mathbb{D}}).$$
More exactly, the bijection is as follows:
\begin{equation*}
\mathscr{D}(\Sigma_{\mathbb{D}}) \to \mathscr{D}(\Delta_{\mathbb{D}}):\quad
f(x_1,x_2,\ldots,x_n) \mapsto f(e_1(u), q^2 e_2(u), \ldots, q^{n(n-1)}
e_n(u)).
\end{equation*}

Thereby,
\begin{equation}\label{inv_iso_nc}
\mathscr{D}(\Delta_{\mathbb{D}}) \cong
 \mathscr{D}(\Sigma_{\mathbb{D}}) \cong \mathscr{D}(\mathbb{D})^{U_q
\mathfrak{s} (\mathfrak{gl}_n \times \mathfrak{gl}_n)}_q. \end{equation} In the
sequel we do not distinguish between $\mathscr{D}(\Delta_{\mathbb{D}})$ and
$\mathscr{D}(\mathbb{D})^{U_q \mathfrak{s} (\mathfrak{gl}_n \times
\mathfrak{gl}_n)}_q$.


Recall the definition of a multiple Jackson integral with 'base' $q^{-2}$ (see
\cite{StokmanSIAM}):
\begin{equation}\label{Jackson_minus}
\int_{q^{-2(n-1)}}^\infty \int_{q^{-2(n-2)}}^{q^2u_n} \ldots
\int_{1}^{q^2u_2} \phi(u) d_{q^{-2}}u_1 \ldots d_{q^{-2}}u_n
\stackrel{\mathrm{def}}{=} (1-q^2)^n \sum\limits_{\lambda \in \Lambda_n}
 \phi(q^{-2(\lambda+\delta)}) q^{-2|\lambda+\delta|}.
\end{equation}

\begin{proposition}\label{intourns_prop}
The restriction of the invariant integral (\ref{invin_d}) to the space
$\mathscr{D}(\mathbb{D})^{U_q \mathfrak{s} (\mathfrak{gl}_n \times
\mathfrak{gl}_n)}_q$ is $$\int f(x_1,x_2,\ldots,x_n)\ d\nu =$$
\begin{equation*} \label{intourns}
 \mathcal{N}  \int_{q^{-2(n-1)}}^\infty
\int_{q^{-2(n-2)}}^{q^2u_n} \ldots \int_{1}^{q^2u_2} f(e_1(u), q^2 e_2(u),
\ldots, q^{n(n-1)} e_n(u)) \Delta(u)^2 \; d_{q^{-2}}u_1 d_{q^{-2}}u_2 \ldots
d_{q^{-2}}u_n,
\end{equation*}
where $\Delta(u) = \prod\limits_{1\leq i<j\leq n}(u_i-u_j)$, $\mathcal{N} =
(1-q^2)^{n(n-1)} q^{n(n-1)} \Delta(q^{-2\delta})^{-2}$ is a positive
constant.
\end{proposition}

The constant $\mathcal{N}$ can be found easily by calculating the integral
for the element $f_0$:
$$
\int f_0 \; d\nu = (1-q^2)^{n^2} = \mathcal{N} (1-q^2)^n
\Delta(q^{-2\delta})^2 q^{-2|\delta|}.
$$
So,
\begin{equation}\label{N_norm}
\mathcal{N} = (1-q^2)^{n(n-1)} q^{n(n-1)} \Delta(q^{-2\delta})^{-2}.
\end{equation}

{\bf Proof.} 
Consider the integral
\begin{equation*}
\widetilde{\eta}: f \mapsto \int_{q^{-2(n-1)}}^\infty
\int_{q^{-2(n-2)}}^{q^2u_n} \ldots \int_{1}^{q^2u_2} f(e_1(u), q^2 e_2(u),
\ldots, q^{n(n-1)} e_n(u)) \Delta(u)^2 \; d_{q^{-2}}u_1 d_{q^{-2}}u_2 \ldots
 d_{q^{-2}}u_n.
\end{equation*}


Let us show that integrals $\eta$ and $\widetilde{\eta}$ are equal up to a
multiplicative constant on the space $\mathscr{D}(\mathbb{D})^{U_q
\mathfrak{s} (\mathfrak{gl}_n \times \mathfrak{gl}_n)}_q$ (the normalizing
constant is calculated in \eqref{N_norm}).

Let us compute $\widetilde{\eta}(f)$:
$$\widetilde{\eta}(f) =\mathrm{const} \int_{q^{-2(n-1)}}^\infty  \ldots
\int_{1}^{q^2u_2} f(e_1(u), q^2 e_2(u), \ldots, q^{n(n-1)} e_n(u))
\Delta(u)^2 \; d_{q^{-2}}u_1 \ldots d_{q^{-2}}u_n=$$

$$= \mathrm{const}
\sum\limits_{\lambda \in \Lambda_n} f(e_1(q^{-2(\lambda +\delta)}), q^2
e_2(q^{-2(\lambda +\delta)}), \ldots, q^{n(n-1)} e_n(q^{-2(\lambda
+\delta)})) \; \Delta(q^{-2(\lambda + \delta)})^2 \; q^{-2|\lambda +
\delta|}.$$

Let us also compute $\eta(f)$:
$$\eta(f) = \mathrm{const} \; \mathrm{tr}(T_F(f)q^{-2\check{\rho}}) =
 \mathrm{const} \sum\limits_{\lambda \in \Lambda_n}
\mathrm{tr} \;(\; T_F(f)|_{\mathcal{H}_\lambda} \;
q^{-2\check{\rho}}|_{\mathcal{H}_\lambda} \;) = $$ $$= \mathrm{const}
\sum\limits_{\lambda \in \Lambda_n} \; d_{\lambda} \; f(e_1(q^{-2(\lambda
+\delta)}), q^2 e_2(q^{-2(\lambda +\delta)}), \ldots, q^{n(n-1)}
e_n(q^{-2(\lambda +\delta)})),
$$ where $d_{\lambda} = \mathrm{tr} \;(\;
q^{-2\check{\rho}}|_{\mathcal{H}_\lambda} \; )$. In the last computation we
essentially use the fact that the operators $T_F(f)$, $f \in
\mathscr{D}(\mathbb{D})^{U_q \mathfrak{s} (\mathfrak{gl}_n \times
\mathfrak{gl}_n)}_q$ are scalar on each $\mathcal{H}_\lambda$.

Introduce the notation $$H_0 = n H_n + \sum\limits_{j=1}^{n-1} j H_j +
\sum\limits_{j=1}^{n-1} j H_{2n-j},$$ then
$$-2\check{\rho} = -n H_0 - \sum\limits_{j=1}^{n-1} j(n-j) H_j -
\sum\limits_{j=1}^{n-1} j(n-j) H_{2n-j}.$$

Consider the subalgebra in $U_q\mathfrak{s}(\mathfrak{gl}_n\times
\mathfrak{gl}_n)$ generated by $\{E_j,F_j,K^{\pm1}_j\}_{j \neq n}$. It is
isomorphic to $U_q\mathfrak{sl}_n\otimes U_q \mathfrak{sl}_n$. The
restriction of the representation of $U_q \mathfrak{s} (\mathfrak{gl}_n
\times \mathfrak{gl}_n)$ in $\mathcal{H}_\lambda$ to the subalgebra $U_q
\mathfrak{sl}_n\otimes U_q \mathfrak{sl}_n$ is equivalent to the
representation $\pi\boxtimes\pi$, where $\pi$ is an irreducible
representation of $U_q \mathfrak{sl}_n$ with highest weight
$(\lambda_1-\lambda_2,\lambda_2-\lambda_3,\ldots,\lambda_{n-1}-\lambda_n)$.
Consequently (see \cite[\S 7.1.4]{KlSch}),
$$d_{\lambda} =
\mathrm{tr}\; ( \; q^{-nH_0}|_{\mathcal{H}_\lambda}) \; ( \; \mathrm{tr}(
\pi(q^{-2 \check{\rho}^{(n)}})|_{\mathcal{H}^{(n)}_\lambda} ) \; )^2 =
q^{-2|\lambda|} S_\lambda(q^{-2\delta})^2,$$ where $\check{\rho}^{(n)} =
\sum\limits_{j=1}^{n-1} j(n-j) H_j$, and
$$S_\lambda(z_1,z_2,\ldots,z_n) = \frac
{ \mathrm{det}(z_i^{\lambda_j + j-1})_{i,j=1,2,\ldots,n}
} { \mathrm{det}(z_i^{j-1})_{i,j=1,2,\ldots,n}}
$$ is the Schur polynomial \cite[\S 1.3]{Mac95}. So $S_\lambda(q^{-2\delta}) =
\frac{\Delta(q^{-2(\lambda + \delta)})}{\Delta(q^{-2\delta})}$, and $$
\eta(f) = \mathrm{const} \sum\limits_{\lambda \in \Lambda_n} \; d_{\lambda}
\; f(e_1(q^{-2(\lambda +\delta)}), q^2 e_2(q^{-2(\lambda +\delta)}), \ldots,
q^{n(n-1)} e_n(q^{-2(\lambda +\delta)})) = $$ $$= \mathrm{const} \;
\sum\limits_{\lambda \in \Lambda_n} \; q^{-2|\lambda|} \;
\Delta(q^{-2(\lambda + \delta)})^2 \; f(e_1(q^{-2(\lambda +\delta)}), q^2
e_2(q^{-2(\lambda +\delta)}), \ldots, q^{n(n-1)} e_n(q^{-2(\lambda
+\delta)})).$$

Now it is obvious that the integrals $\eta$ and $\widetilde{\eta}$ are equal up to
a multiplier. \hfill $\Box$

\section{Spherical functions on quantum Grassmanian}

Consider the involution $\star$ in $U_q \mathfrak{sl}_{2n}$ determined by
\begin{equation*}
(K_j^{\pm 1})^\star=K_j^{\pm 1},\qquad E_j^\star=K_jF_j,\qquad
F_j^\star=E_jK_j^{-1}.
\end{equation*}
Then $U_q\mathfrak{su}_{2n} = (U_q\mathfrak{sl}_{2n}, \star)$ is an $*$-Hopf
algebra. It is a quantum analogue of $U \mathfrak{su}_{2n}
\otimes_{\mathbb{R}} \mathbb{C}$.

Consider also the involution $\star$ in $\mathbb{C}[SL_{2n}]_q$ determined
by
$$
t_{ij}^\star = (-q)^{j-i} t^{\wedge 2n-1}_{ \{1,2,\ldots,2n\} \backslash
\{i\}, \{1,2,\ldots,2n\} \backslash \{j\}}.
$$
The $*$-Hopf algebra $\mathbb{C}[SU_{2n}]_q \stackrel{\mathrm{def}}{=}
(\mathbb{C}[SL_{2n}]_q,\star)$ is a $U_q\mathfrak{su}_{2n}$-module $*$-Hopf
algebra. It is a well-known quantum analogue of the algebra of regular
functions on the Lie group $SU_{2n}$ (see \cite{Woron1}, \cite{Woron91}).

It is well known that in the classical case the Cartan duality between
compact and noncompact Hermitian symmetric spaces allows one to predict some
results of harmonic analysis in the noncompact case using the easier compact
case. In this \itemii $\ $ we explore this observation. We construct a
family of difference operators for the quantum Grassmanians. These operators
are obtained using the action of the center of $U_q\mathfrak{sl}_{2n}$.
Afterwards, our construction allows us to introduce difference operators in
the case of quantum matrix ball.

\subsection{Spherical functions}

It is well known that for any finite-dimensional irreducible $U_q
\mathfrak{sl}_{2n}$-module $V$
$$\dim V^{U_q\mathfrak{s}(\mathfrak{gl}_n \times \mathfrak{gl}_n)}\leq1.$$
Hence $(U_q \mathfrak{sl}_{2n}, U_q \mathfrak{s} (\mathfrak{gl}_n \times
\mathfrak{gl}_n))$ is a "quantum Gelfand pair". As in the classical case, let us
define a simple finite-dimensional weight $U_q \mathfrak{sl}_{2n}$-module to
be spherical, if $\dim V^{U_q\mathfrak{s}(\mathfrak{gl}_n \times
\mathfrak{gl}_n)}=1.$

\begin{remark}
It is well known (\cite{StDij}, \cite[theorem 4.4.1]{StokmanPhD}) that a
simple finite-dimensional weight $U_q \mathfrak{sl}_{2n}$-module is $U_q
\mathfrak{s}(\mathfrak{gl}_n \times \mathfrak{gl}_n)$-spherical if and only
if its highest weight has the following form:
\begin{equation*}
\widehat{\lambda} = (\lambda_1 - \lambda_2, \lambda_2 - \lambda_3,
\ldots,\lambda_{n-1} - \lambda_n, 2\lambda_n, \lambda_{n-1} - \lambda_n,
\ldots, \lambda_2 - \lambda_3, \lambda_1 - \lambda_2),\qquad \lambda \in
\Lambda_n. \end{equation*} We will denote by $L_\lambda$ the $U_q
\mathfrak{sl}_{2n}$-module with highest weight $\widehat{\lambda}$.
\end{remark}

A scalar product\footnote{A sesquilinear positive definite Hermitian
symmetric form.} $(\cdot,\cdot)$ in $V$ is called $U_q
\mathfrak{su}_{2n}$-invariant if for any $\xi\in U_q \mathfrak{sl}_{2n}$ and
for any $v_1,v_2\in V$ $$(\xi v_1,v_2)=(v_1,\xi^\star v_2).$$ Any spherical
$U_q \mathfrak{sl}_{2n}$-module $V$ can be equipped with a $U_q
\mathfrak{su}_{2n}$-invariant scalar product. Fix ${v \in V^{U_q
\mathfrak{s} (\mathfrak{gl}_n \times \mathfrak{gl}_n)}}$ by the requirement
$(v,v)=1$. Recall (see \cite[\S 11.6.4]{KlSch}) that the matrix element
$\varphi_V(\xi)=(\xi v, v)$ corresponding to $U_q \mathfrak{s}
(\mathfrak{gl}_n \times \mathfrak{gl}_n)$-invariant vector is called the
spherical function on the quantum group $SU_{2n}$ corresponding to $V$.

Thus $\varphi_V$ is a $U_q \mathfrak{s} (\mathfrak{gl}_n \times
\mathfrak{gl}_n)$-biinvariant element of $\mathbb{C}[SU_{2n}]_q$ such that
\begin{equation*}
\varphi_V(1) = 1.
\end{equation*}
The lemma below follows from the results of \cite{KorSoib}.
\begin{lemma}
\begin{equation*}
(\varphi_V)^\star=\varphi_V.
\end{equation*}
\end{lemma}

\bigskip
It follows from Proposition 7 of \cite{BKV} and Lemma 1 of \cite{BerS} that
the subalgebra of $U_q \mathfrak{s} (\mathfrak{gl}_n \times
\mathfrak{gl}_n)$-biinvariant functions in $\mathbb{C}[SU_{2n}]_q$ is
generated by the pairwise commuting elements $x_1,x_2,\ldots,x_n$. In
particular, every spherical function $\varphi_V$ is a polynomial in
$x_1,x_2,\ldots,x_n$. Denote by $\varphi_\lambda(x_1,x_2,\ldots,x_n)$ the
spherical function corresponding to the module $L_\lambda$. In this \itemi$
$ we will find an exact formula for $\varphi_\lambda(x_1,x_2,\ldots,x_n)$.

\subsubsection{Little $q$-Jacobi polynomials}

We will use the following partial order on $\Lambda_n$
\begin{equation*}
\eta\leq\lambda \stackrel{\mathrm{def}}{\Longleftrightarrow}
\sum\limits_{j=1}^k \eta_j \leq \sum\limits_{j=1}^k \lambda_j, \qquad k =
1,2, \ldots, n.
\end{equation*}
As usual, $ \eta<\lambda \stackrel{\mathrm{def}}{\Longleftrightarrow} \eta
\leq \lambda \quad \& \quad \eta \neq \lambda$.

Introduce the short notation $\mathbf{1}^k =
(\underbrace{1,\ldots,1}_k,0,\ldots,0)$. Let us denote by $m_\lambda$ the
monic symmetric polynomial
$$
m_\lambda(z_1,z_2,\ldots,z_n) = \sum\limits_{w \in S_n} z_{w(1)}^{\lambda_1}
z_{w(2)}^{\lambda_2} \ldots z_{w(n)}^{\lambda_n}.$$

Let $P_\lambda$ be a unique symmetric polynomial which satisfies the
following two conditions

\begin{equation*}
\begin{tabular}{c}
$1) \qquad P_\lambda(z)=m_\lambda(z)+\sum\limits_{\eta < \lambda}
d_{\lambda,\eta}m_\eta(z), \quad d_{\lambda,\eta} \in \mathbb{R},$\\$ 2)
\qquad \int_{0}^{q^2} \ldots \int_{0}^{q^2z_2} P_\lambda(z) m_\eta(z)
\Delta(z)^2 d_{q^2}z_1 \ldots d_{q^2}z_n = 0, \quad \eta<\lambda,$
\end{tabular}
\end{equation*}
where the multiple Jackson integral (cf. \eqref{Jackson_minus}) is defined as
$$
\int_{0}^{q^2} \ldots \int_{0}^{q^2z_2} \phi(z) d_{q^2}z_1 \ldots d_{q^2}z_n
= (1-q^2)^n \sum\limits_{\lambda \in \Lambda_n}
 \phi(q^{2(\lambda+\delta+ \mathbf{1}^n)}) q^{2|\lambda+\delta+
 \mathbf{1}^n|}.
$$

\begin{remark}
It is easy to see that $$P_\lambda(z) = P_\lambda(z; 0,0; q^2),$$ where
$P_\lambda(z; a,b; q)$ are Little $q$-Jacobi polynomials (see
\cite{StokmanSIAM}).
\end{remark}

Let $\widetilde{P}_\lambda$ be a polynomial such that
$$P_\lambda(z) =  \widetilde{P}_\lambda(e_1(z), q^2 e_2(z), \ldots, q^{n(n-1)}
e_n(z)).$$ From the results of \itemii $ $ \ref{x_spectr_integral} and
\cite{StDij}, \cite[Theorem 4.7.5]{StokmanPhD} one can deduce the following
theorem:
\begin{theorem}
The spherical function $\varphi_\lambda$ is equal (up to a multiplicative
constant) to
\begin{equation*}
\widetilde{P}_\lambda(x_1, x_2, \ldots, x_n).
\end{equation*}
\end{theorem}

Denote the fundamental spherical weights by
\begin{equation*}
\mu_k = \widehat{\mathbf{1}^k}, \qquad k \in \{1,2,\ldots,n\},
\end{equation*}
and denote by
\begin{equation*}P^{\mathrm{spher}}_+ = \bigoplus\limits_{k=1}^n
\mathbb{Z}_+ \mu_k = \{ \widehat{\lambda}\;|\; \lambda \in \Lambda_n \}
\end{equation*} the set of positive spherical weights, and by
\begin{equation}\label{p_lambda}
P^{\mathrm{spher}} = \bigoplus\limits_{k=1}^n \mathbb{Z} \mu_k = \{
\widehat{\lambda} \;|\; \lambda \in \mathbb{Z}^n \}
\end{equation} the set of all
spherical weights.

Stokman proved the following formula in \cite[Proposition 5.9]{StokmanSIAM}:
$$ P_\lambda(z_1, z_2, \ldots, z_n) = \Delta(z)^{-1} \sum\limits_{w \in S_n}
\mathrm{sign}(w) \prod\limits_{i=1}^n P_{(\lambda+\delta)_{w(i)}}(z_i),
$$
where $P_m(z)$ are Little $q$-Jacobi polynomials in one variable.

Recall the 'coordinates' $u_1,u_2,\ldots,u_n$ appeared in
\eqref{x_e_link_inj}.
\begin{corollary}\label{det_formula} Let $\lambda \in
\Lambda_n$. Then
$$\varphi_\lambda(u) =  \mathrm{const} \; P_\lambda(u) =
\mathrm{const} \Delta(u)^{-1} \sum\limits_{w \in S_n} \mathrm{sign}(w)
\prod\limits_{i=1}^n P_{d(\lambda,w,i)}(u_i),$$ where
$d(\lambda,w,i)=(\lambda+\delta)_{w(i)}\in \mathbb{Z}$.
\end{corollary}


\subsection{Difference operators and the action of the center of
$U_q^{\mathrm{ext}} \mathfrak{sl}_{2n}$ }\label{difoper_compact}

Let $a_{ij}$ be the Cartan matrix of the Lie algebra $\mathfrak{sl}_{2n}$. Denote
by $\alpha_i$, $i = 1,2,\ldots, 2n-1$, simple roots such that
$\alpha_i(H_j)=a_{ji}$ and by $\Phi$ the root system of the Lie algebra
$\mathfrak{sl}_{2n}$.

In this \itemii$\ $ we will consider the action of the center of
$U_q\mathfrak{sl}_{2n}$ in weight modules. Note that it is more convenient
to use the center of the extended quantum universal enveloping algebra
$U_q^{\mathrm{ext}} \mathfrak{sl}_{2n}$. Essentially,
$U_q^{\mathrm{ext}}\mathfrak{sl}_{2n}$ can be obtained from
$U_q\mathfrak{sl}_{2n}$ by adding the elements
\begin{equation*}
K_\lambda=K_1^{a_1}K_2^{a_2}\ldots K_{2n-1}^{a_{2n-1}}, \qquad
\mathbf{\lambda}=\sum_{i=1}^{2n-1} a_i\alpha_i
\end{equation*}
for all $\lambda$ in the weight lattice $P$. In particular, the action of
$U_q\mathfrak{sl}_{2n}$ in any weight module admits a unique extension to the
action of $U_q^{\mathrm{ext}}\mathfrak{sl}_{2n}$. Denote by
$Z(U_q^{\mathrm{ext}} \mathfrak{sl}_{2n})$ the center of the extended universal
enveloping algebra.

Recall some definitions, cf. \cite{Bou4-6Eng}. Consider the real linear span
$\mathfrak{h}^*_\mathbb{R}$ of all simple roots of the Lie algebra
$\mathfrak{sl}_{2n}$. It is well known that there is a positive definite scalar
product $(\cdot, \cdot)$ in $\mathfrak{h}^*_\mathbb{R}$. Denote by
$(\mathfrak{h}^*_\mathbb{R})^- \subset \mathfrak{h}^*_\mathbb{R}$ the real
subspace spanned by the strictly orthogonal noncompact positive roots
$$\gamma_k = \alpha_k + \alpha_{k+1} + \ldots + \alpha_{2n-k-1} +
\alpha_{2n-k}, \qquad k \in \{ 1, 2, \ldots, n \},$$ and by
$(\mathfrak{h}^*_\mathbb{R})^+ \subset \mathfrak{h}^*_\mathbb{R}$ its
orthogonal complement. It is well known that the orthogonal projection of the
root system $\Phi$ to $(\mathfrak{h}^*_\mathbb{R})^-$ is a root system of type
$C_n$ and it is called the system of restricted roots $\Phi^{\mathrm{res}}$. The
Weyl group $W^{\mathrm{res}}$ of the root system $\Phi^{\mathrm{res}}$ is
called the restricted Weyl group.

Let $\mathbb{C}[ P^{\mathrm{spher}}]_q$ be an algebra generated by the
following functions on $P^{\mathrm{spher}}$:
\begin{equation*} \lambda \mapsto q^{(\eta, \lambda)}, \quad \eta \in
P^{\mathrm{spher}}.
\end{equation*}
This algebra is naturally isomorphic to the group algebra of the lattice
$P^{\mathrm{spher}}$. Denote by $\mathbb{C}[
P^{\mathrm{spher}}]_q^{W^{\mathrm{res}}}$ the subalgebra of
$W^{\mathrm{res}}$-invariants in $\mathbb{C}[ P^{\mathrm{spher}}]_q$:
\begin{equation*}
\mathbb{C}[ P^{\mathrm{spher}}]_q^{W^{\mathrm{res}}} = \{ f \in \mathbb{C}[
P^{\mathrm{spher}}]_q \quad |
 \quad f(w \lambda) = f(\lambda)\quad \text{for all} \quad w \in
W^{\mathrm{res}}, \quad \lambda \in P^{\mathrm{spher}}\}.
\end{equation*}

Here we provide a well-known description of the image of the center
$Z(U_q^{\mathrm{ext}} \mathfrak{sl}_{2n})$ under the Harish-Chandra
homomorphism $ \gamma^{\mathrm{spher}} : Z(U_q^{\mathrm{ext}}
\mathfrak{sl}_{2n}) \to \mathbb{C}[P^{\mathrm{spher}}]_q$ (see
\cite{Baldoni_Frajria}).

\begin{proposition}\label{general_center_image}
The image of $Z(U_q^{\mathrm{ext}} \mathfrak{sl}_{2n})$ under the
Harish-Chandra homomorphism is the subalgebra $\mathbb{C}[
P^{\mathrm{spher}}]_q^{W^{\mathrm{res}}}$.
\end{proposition}

Set for $\lambda \in \mathbb{C}^n$ $$a(\lambda+\delta)
\stackrel{\mathrm{def}}{=} (a(\lambda_1+n-1), a(\lambda_2+n-2), \ldots,
a(\lambda_n)),$$ where
\begin{equation}\label{eigenvalue}
a(l) = \frac{(1-q^{-2l})(1-q^{2l+2})} {(1-q^2)^2}, \qquad l \in \mathbb{C}.
\end{equation}

\begin{proposition}\label{C_phi}
There are such elements $C_k \in Z(U_q^{\mathrm{ext}} \mathfrak{sl}_{2n})$,
$k=1,2,\ldots, n$, that
\begin{equation}\label{C_k_def}
C_k \varphi_\lambda = e_k(a(\lambda+\delta)) \varphi_\lambda,\qquad \lambda
\in \Lambda_n.
\end{equation}
\end{proposition}
{\bf Proof.} Consider the mapping:
$$\Lambda_n \to \mathbb{R}^n,\qquad \lambda \mapsto \eta(\lambda) = (\lambda_1-\frac{2n-1}{2},\lambda_2-\frac{2n-3}{2},
\ldots,\lambda_{n-1}-\frac{3}{2}, \lambda_n - \frac{1}{2} ) .$$ Then
$$\widehat{\eta(\lambda)} = \widehat{\lambda} - \rho,$$
where
\begin{equation*}
\widehat{\lambda} = (\lambda_1 - \lambda_2, \lambda_2 - \lambda_3,
\ldots,\lambda_{n-1} - \lambda_n, 2\lambda_n, \lambda_{n-1} - \lambda_n,
\ldots, \lambda_2 - \lambda_3, \lambda_1 - \lambda_2) \in P.
\end{equation*}

We need the following functions on $P^{\mathrm{spher}}$:
\begin{equation}\label{rho_shifted_spectr}
\psi_k: \widehat{\lambda} \mapsto e_k(a(\eta(\lambda)+\delta)),\qquad k \in
\{1,2,\ldots,n\},
\end{equation}
$\lambda\in \mathbb{Z}^n$ is uniquely defined by the spherical weight
$\widehat{\lambda} \in P^{\mathrm{spher}}$, see \eqref{p_lambda}.

Due to Proposition \ref{general_center_image} we only need to check the
$W^{\mathrm{res}}$-invariance of the functions $\psi_k$.

It is easy to see that $$e_k(a(\eta(\lambda)+\delta)) = $$ $$(1-q^2)^{-2k}
 e_k((1-q^{-2\lambda_1+1})(1-q^{2\lambda_1+1}),
(1-q^{-2\lambda_2+1})(1-q^{2\lambda_2+1}), \ldots,
(1-q^{-2\lambda_n+1})(1-q^{2\lambda_n+1})).$$ Besides,
$$\widehat{\lambda} = \lambda_1 \gamma_1 + \lambda_2 \gamma_2 +
\ldots + \lambda_n \gamma_n.$$ As the group $W^{\mathrm{res}}$ acts on
$\gamma_k$ by permutations and sign changes, the function
\eqref{rho_shifted_spectr} is $W^{\mathrm{res}}$-invariant. \hfill $\Box$

\bigskip

Let $\mathcal{L}_k$ be the linear operator in $\mathbb{C}[SL_{2n}]_q$ defined
by $\mathcal{L}_k f = C_k f$.

The action of $U_q^{\mathrm{ext}}\mathfrak{sl}_{2n}$ in the space of
$U_q\mathfrak{s}(\mathfrak{gl}_n\times \mathfrak{gl}_n)$-biinvariant functions
determines the homomorphism
\begin{equation*}
Z(U_q^{\mathrm{ext}}\mathfrak{sl}_{2n})\to
\mathrm{End}(\mathbb{C}[u_1,u_2,\ldots,u_n]^{S_n}),
\end{equation*}
as
\begin{equation}\label{inv_iso}\mathbb{C}[u_1,u_2,\ldots,u_n]^{S_n} \cong
\mathbb{C}[x_1,x_2,\ldots,x_n] \cong \mathbb{C}[SL_{2n}]^{(U_q
\mathfrak{s}(\mathfrak{gl}_n \times \mathfrak{gl}_n))^{\mathrm{op}} \otimes
U_q \mathfrak{s}(\mathfrak{gl}_n \times \mathfrak{gl}_n)}_q,\end{equation}
see \itemii $ $ \ref{x_spectr_integral}. Here we will describe the action of
the linear operators $\mathcal{L}_1$, $\mathcal{L}_2$, ..., $\mathcal{L}_n$
in the space \eqref{inv_iso}.


Let us define the difference operator $\Box_{u_i}$ in the space
$\mathbb{C}[u_1,u_2,\ldots,u_n]$ with
\begin{equation}\label{sqr_box} \Box_{u_i} f(u_1, \ldots, u_n) = D_{u_i} u_i (1-q^{-1} u_i)
D_{u_i} f(u_1, \ldots, u_n),\end{equation} where $D_{u_i} f(u_1, \ldots,
u_n) = \frac{f(u_1, \ldots, u_{i-1}, q^{-1}u_i, u_{i+1}, \ldots, u_n)-f(u_1,
\ldots, u_{i-1}, qu_i, u_{i+1}, \ldots, u_n)}{q^{-1}u_i-qu_i}$.

\begin{proposition}
\begin{equation}\label{diff_oper}
\mathcal{L}_k|_{\mathbb{C}[u_1, u_2, \ldots, u_n]^{S_n}} \; =
\frac{1}{\Delta(u)} \; e_k(\Box_{u_1}, \ldots, \Box_{u_n}) \; \Delta(u).
\end{equation}
\end{proposition}
{\bf Proof.} In \itemii $\ $ \ref{Plansh_1} it will be showed that in the
case of one variable $$\Box_u \; \varphi_l(u) = a(l) \varphi_l(u).
$$

From \eqref{eigenvalue} and the determinant decomposition described in
Corollary \ref{det_formula} it follows that
$$
\frac{1}{\Delta(u)} \; e_k(\Box_{u_1}, \ldots, \Box_{u_n}) \; \Delta(u) \;
\varphi_\lambda(u) = e_k(a(\lambda+\delta)) \varphi_\lambda(u),\qquad
\lambda \in \Lambda_n.
$$

The equality \eqref{diff_oper} follows from Proposition \ref{C_phi}, as the set
$\{\varphi_\lambda\}_{\lambda \in \Lambda_n}$ is a basis of the vector space
$\mathbb{C}[SL_{2n}]^{(U_q \mathfrak{s}(\mathfrak{gl}_n \times
\mathfrak{gl}_n))^{\mathrm{op}} \otimes U_q \mathfrak{s}(\mathfrak{gl}_n
\times \mathfrak{gl}_n)}_q$. \hfill $\Box$

\section{Plancherel measure for the quantum matrix ball}
\subsection{The Plancherel measure for a family of operators
$\mathcal{L}^{\mathrm{radial}}_1$, $\mathcal{L}^{\mathrm{radial}}_2$, ...,
$\mathcal{L}^{\mathrm{radial}}_n$}\label{double_inv_plansh}

\subsubsection{Linear operators $\mathcal{L}^{\mathrm{radial}}_1$,
$\mathcal{L}^{\mathrm{radial}}_2$, ..., $\mathcal{L}^{\mathrm{radial}}_n$ in the
space $L^2(\Delta_{\mathbb{D}},d\nu_q)$}\label{center_nc}

Let us consider the elements $$C_1, C_2, \ldots,C_n \in Z(U_q^{\mathrm{ext}}
\mathfrak{sl}_{2n})$$ defined in \eqref{C_k_def}. Let also $\mathcal{L}_k$ be the
linear operator in $\mathscr{D}(\mathbb{D})_q$ defined by
$$
\mathcal{L}_k f = C_k f.
$$

Now we describe the restriction of the linear operator $\mathcal{L}_k$, $k=1,2,
\ldots, n$, to the space $\mathscr{D}(\mathbb{D})^{U_q \mathfrak{s}
(\mathfrak{gl}_n \times \mathfrak{gl}_n)}_q$ of $U_q \mathfrak{s}
(\mathfrak{gl}_n \times \mathfrak{gl}_n)$-invariants in
$\mathscr{D}(\mathbb{D})_q$.

Let us introduce the short notation $\mathcal{L}_k^{\mathrm{radial}}$ for the
restriction of $\mathcal{L}_k$ to $\mathscr{D}(\mathbb{D})^{U_q \mathfrak{s}
(\mathfrak{gl}_n \times \mathfrak{gl}_n)}_q$.

\begin{proposition}\label{qlimittrans}
\begin{equation}\label{diff_oper_nc}
\mathcal{L}_k^{\mathrm{radial}} \; = \frac{1}{\Delta(u)} \; e_k(\Box_{u_1},
\ldots, \Box_{u_n}) \; \Delta(u),
\end{equation}
where $\Box_{u_j}$ are the difference operators in the vector space
\eqref{inv_iso_nc} defined by the same formula as in \eqref{sqr_box}.
\end{proposition}

{\bf Proof.} Following \itemii $ $ \ref{x_spectr_integral}, the vector space
of $U_q \mathfrak{s} (\mathfrak{gl}_n \times \mathfrak{gl}_n)$-invariants in
$\mathscr{D}(\mathbb{D})_q$ can be identified with the space
$\mathscr{D}(\Sigma_{\mathbb{D}})$ of functions on $\Sigma_{\mathbb{D}}$
with finite support. Using Lemma \ref{spectr_bij} one can obtain that the
vector space of $U_q \mathfrak{s} (\mathfrak{gl}_n \times
\mathfrak{gl}_n)$-invariants in $\mathscr{D}(\mathbb{D})_q$ is canonically
isomorphic to the space $\mathscr{D}(\Delta_{\mathbb{D}})$ of functions on
$\Delta_{\mathbb{D}}$ with finite support. Consider the point-wise
convergence topology on $\Delta_{\mathbb{D}}$.

The space of symmetric polynomials in $\Delta_{\mathbb{D}}$ is dense in
topological space $\mathcal{F}(\Delta_{\mathbb{D}})$ of functions on
$\Delta_{\mathbb{D}}$, and the equation \eqref{diff_oper_nc} takes place for
symmetric polynomials following \eqref{diff_oper} and \eqref{inv_iso}.

The linear operators in both parts of equation \eqref{diff_oper_nc} can be
extended continuously from the space of symmetric polynomials in
$\Delta_{\mathbb{D}}$ to $\mathcal{F}(\Delta_{\mathbb{D}})$, and the equation
\eqref{diff_oper_nc} takes place for the whole space
$\mathcal{F}(\Delta_{\mathbb{D}})$. \hfill $\Box$

Now we recall the measure on $\Delta_{\mathbb{D}}$:
\begin{equation}\label{measure_n}
d\nu_q(u)= \mathcal{N} \; \Delta(u)^2 d_{q^{-2}}u_1 d_{q^{-2}}u_2 \ldots
d_{q^{-2}}u_n,
\end{equation}
where $\mathcal{N}$ is defined in (\ref{N_norm}). It is the restriction of
the invariant measure to the space $\mathscr{D}(\mathbb{D})^{U_q
\mathfrak{s} (\mathfrak{gl}_n \times \mathfrak{gl}_n)}_q$, which we already
identified with $\mathscr{D}(\Delta_{\mathbb{D}})$, see Proposition
\ref{intourns_prop}.

Let us introduce the Hilbert space $L^2(\Delta_{\mathbb{D}},d\nu_q)$ of
functions on the set $\Delta_{\mathbb{D}}$ which satisfy
$$
\int\limits_{\Delta_{\mathbb{D}}} |f(u)|^2 d\nu_q(u) <\infty,
$$
where
$$
(f,g)=\int\limits_{\Delta_{\mathbb{D}}} \overline{g(u)} f(u) d\nu_q(u).
$$

It will be proved in the sequel (Lemma \ref{oper}) that the linear operators
$\mathcal{L}^{\mathrm{radial}}_1$, $\mathcal{L}^{\mathrm{radial}}_2$, ...,
$\mathcal{L}^{\mathrm{radial}}_n$ can be continuously extended to the
bounded pairwise commuting self-adjoint operators in
$L^2(\Delta_{\mathbb{D}},d\nu_q)$.

Our goal is to find a Plancherel measure $d\Sigma$ on the joint spectrum of
commuting self-adjoint linear operators $\mathcal{L}^{\mathrm{radial}}_1$,
$\mathcal{L}^{\mathrm{radial}}_2$, ..., $\mathcal{L}^{\mathrm{radial}}_n$ and a
unitary operator $\mathscr{F}: L^2(\Delta_{\mathbb{D}},d\nu_q) \to
L^2(d\Sigma)$ which provides a unitary equivalence between the operators
$\mathcal{L}^{\mathrm{radial}}_1$, $\mathcal{L}^{\mathrm{radial}}_2$, ...,
$\mathcal{L}^{\mathrm{radial}}_n$ and the operators of multiplication by
independent variable, such as $\mathscr{F} f_0 = 1$.

The element $f_0 \in L^2(\Delta_{\mathbb{D}},d\nu_q)$ is a cyclic vector
under the action of $\mathcal{L}^{\mathrm{radial}}_1$,
$\mathcal{L}^{\mathrm{radial}}_2$, ..., $\mathcal{L}^{\mathrm{radial}}_n$
(one can prove it explicitly, see \itemii$\ $\ref{f0_cyc}). However, it
follows from the isometry of the operator $\mathscr{F}$ and
Remarks \ref{fo_im} and \ref{independ_var}.

The considered problems are typical for the theory of commutative operator
$*$-algebras with a cyclic vector \cite[p. 570-571]{Naimark}, \cite[p.
103]{Takesaki}.

\subsubsection{The cyclic vector $f_0$}\label{f0_cyc}

Here we discuss the fact that 
the element $f_0 \in L^2(\Delta_{\mathbb{D}},d\nu_q)$ is a cyclic vector
under the action of the operators $\mathcal{L}^{\mathrm{radial}}_1$,
$\mathcal{L}^{\mathrm{radial}}_2$, ..., $\mathcal{L}^{\mathrm{radial}}_n$.

By direct computation we obtain the following lemma.

\begin{lemma}\label{f_0_r1}
In the case of quantum disk \begin{equation}\label{f0_l4}\Box_u f_0(q^{2k}
u) = c_{k,-1} f_0(q^{2k-2}u) + c_{k,0} f_0 (q^{2k} u) + c_{k,1} f_0
(q^{2k+2} u), \qquad k \in \mathbb{Z}_+,\end{equation} where $c_{k, -1}$,
$c_{k,0}$, $c_{k,1}$ are nonzero constants.
\end{lemma}
For example, $$\square_u f_0(u)=D_u u(1-q^{-1}u)D_u f_0(u)={f_0(u) \over
1-q^2}-{q^2f_0(q^2 u)\over 1-q^2}.$$ Here $f_0(q^{-2} u) = 0$ for $u \in
q^{-2\mathbb{Z}_+}$, so the first term in \eqref{f0_l4} vanishes.

The next lemma follows from the previous one by induction
\begin{lemma}\label{l_act_l5}
$$ \mathcal{L}^{\mathrm{radial}}_i f_0(q^{2(\lambda+\delta)} u ) =
\sum c_\mathbf{d} f_0(q^{2(\lambda+\delta+\mathbf{d})} u),$$
where $\mathbf{d} \in \{ -1, 0, 1 \}^n$, $\mathrm{card} \{j|d_j \neq 0\}
\leq i$ and $c_\mathbf{d} \neq 0$.
\end{lemma}

\begin{lemma}\label{f_0_r1}
The linear span of the action of $\mathcal{L}^{\mathrm{radial}}_1$,
$\mathcal{L}^{\mathrm{radial}}_2$, ..., $\mathcal{L}^{\mathrm{radial}}_n$ on
$f_0$ contains the set of finite functions on
$\mathscr{D}(\Delta_{\mathbb{D}})$.
\end{lemma}
{\bf Sketch of the proof}. Lemma \ref{l_act_l5} implies that the linear span
of the action of $\mathcal{L}^{\mathrm{radial}}_1$,
$\mathcal{L}^{\mathrm{radial}}_2$, ..., $\mathcal{L}^{\mathrm{radial}}_n$ on
$f_0$ contains the set
$$
\mathscr{S}_{\mathbb{D}} = \{ f_0(q^{2(\lambda + \delta)} u)\;|\; \lambda
\in \Lambda_n \}
$$
of characteristic functions of points of $\Delta_{\mathbb{D}}$.

The last lemma implies that $f_0$ is cyclic as the set of finite functions
$\mathscr{D}(\Delta_{\mathbb{D}})$ is dense in $L^2(\Delta_{\mathbb{D}},
d\nu_q)$.

\bigskip

\subsubsection{Example: the quantum disk}\label{Plansh_1}

In this \itemii $ $ we recall the Plancherel measure $d\sigma$ for the quantum
disk found in \cite{VakShkl}.

Consider the Hilbert space $L^2(q^{-2\mathbb{Z}_+})$ of functions on the
geometric series $q^{-2\mathbb{Z}_+}$ which satisfy the condition
$$
\int\limits_1^\infty |f(u)|^2 d_{q^{-2}} u<\infty
$$
with the scalar product
$$
(f,g)=\int\limits_1^\infty \overline{g(u)} f(u) d_{q^{-2}} u.
$$

Recall the notation for the difference operator $\Box_u$, which acts in the space
of functions on geometric series $q^{-2\mathbb{Z}_+}$ by
$$\Box_u f(u)=D_u u(1-q^{-1}u)D_u f(u),$$
where
$$D_u f(u)\mapsto \frac{f(q^{-1}u)-f(qu)}{q^{-1}u-qu}.$$
Then $\mathcal{L}^{\mathrm{radial}}_1 = \Box_u$.

Let us describe the eigenfunctions of the difference operator $\Box_u$.
Introduce the notation
\begin{equation*}\label{phil}
\Phi_l(u)= {_{3}\Phi_{2}}\left(
\begin{array}{c}
q^{-2l},q^{2(l+1)},u\\ q^2,0
\end{array}
;q^2,q^2\right), \qquad l \in \mathbb{C},
\end{equation*}
for a basic hypergeometric function (see \cite{GREng}).

\begin{proposition}\label{own_rk1}(\cite[\S 8]{VakShkl})
$$\Box_u \Phi_l(u)=a(l)\Phi_l(u),$$ where $a(l)$ is defined in
\eqref{eigenvalue}:
\begin{equation*}
 a(l)=\frac{(1-q^{-2l})(1-q^{2l+2})}
{(1-q^2)^2}.
\end{equation*}
\end{proposition}

\begin{remark}
$\Phi_l(1) = 1$.
\end{remark}

\begin{remark}\label{phi_capital}
$\Phi_l(u)$ is equal up to a multiplicative constant to $\varphi_l(u)$.
\end{remark}

Let
$$c(l)=\frac{\Gamma_{q^2}(2l+1)}{(\Gamma_{q^2}(l+1))^2}$$  be a $q$-analogue of the Harish-Chandra $c$-function.
Here
$\Gamma_{q^2}(x)=\frac{(q^2,q^2)_\infty}{(q^{2x},q^2)_\infty}(1-q^2)^{1-x}$
is a well-known $q$-analogue of the Gamma function $\Gamma(x)$.

Let us consider the measure
\begin{equation*}\label{Plansherel_measure}
d \sigma(\rho)= \frac{1}{2 \pi} \cdot \frac{h}{1-q^2} \cdot \frac{d
\rho}{c(-\frac{1}{2} +i \rho) c(-\frac{1}{2} - i \rho)}\end{equation*} on the interval
$[0,\pi/h]$, where $h=-2 \; \mathrm{ln} \; q$.

Consider the operator
$$\mathcal{F}: \; f \mapsto \widehat{f}(\rho)=\int\limits_1^\infty \Phi_{-\frac{1}{2}+i \rho}(u) f(u) d_{q^{-2}} u$$
defined in the space of finite functions on the geometric series
$q^{-2\mathbb{Z}_+}$. It is shown in \cite[Theorem 9.2]{VakShkl} that this
operator can be extended to a unitary operator
$\mathcal{F}:L^2(q^{-2\mathbb{Z}_+}) \to L^2([0,\pi/h],d\sigma)$ such as
$$\mathcal{F} \; \Box_u \; f = a(-\frac{1}{2}+i \rho) \mathcal{F} f,\qquad f \in L^2([0,\pi/h],d\sigma),$$
where $a(l)$ is defined in \eqref{eigenvalue}. The inverse operator
$\mathcal{F}^{-1}$ has the form
$$ \widehat{f}(\rho) \mapsto \int\limits_0^{\pi/h} \widehat{f}(\rho) \Phi_{-\frac{1}{2}+i \rho}(u) d\sigma(\rho).$$

\subsubsection{The quantum matrix ball}\label{Plansh} We will call an
eigenfunction of a difference operator a generalized one if it does not
belong to $L^2$. These functions are used in the sequel for the construction
of the operator $\mathcal{F}$.

Consider the isometric linear operator \footnote{
$L^2(q^{-2\mathbb{Z}^n_+})_q$ is a short notation for
$L^2(\underbrace{q^{-2\mathbb{Z}_+} \times \ldots \times
q^{-2\mathbb{Z}_+}}_{n})_q$ with the product measure, multiplicated by
$\mathcal{N}$. }
$$\mathscr{I}:L^2(\Delta_{\mathbb{D}},d\nu_q) \to
L^2(q^{-2\mathbb{Z}^n_+}),$$
\begin{equation}\label{isooper}
\mathscr{I}: f(u) {\mapsto} \Delta(u) \widetilde{f}(u),
\end{equation}
where $\widetilde{f}$ is defined in the following way: for every $u = (u_1,
\ldots, u_n)$ with $u_i \neq u_j$ for $i \neq j$ there exists a unique
permutation $w \in S_n$ such as $u_{w_1}>u_{w_2}>\ldots>u_{w_n}$. Then
$$ \widetilde{f}(u) = \left \{\begin{array}{cc}
\frac{1 
}{\sqrt{n!}} f(u_{w_1}, \ldots, u_{w_n}),& u_i \neq u_j , \quad i \neq j,\\
0 ,& \text{otherwise}, \end{array}\right.\qquad u_1,u_2,\ldots,u_n \in
q^{-2 \mathbb{Z}_+}.$$

Consider the notation
\begin{equation}\label{raznvyr}\widetilde{\mathcal{L}_k} =
e_k(\Box_{u_1},\Box_{u_2},\ldots,\Box_{u_n}),\qquad
k=1,2,\ldots,n,\end{equation} for the difference operators in
$L^2(q^{-2\mathbb{Z}_+^n} )_q$. Then the following diagram is commutative:
$$  \xymatrix{
L^2(\Delta_{\mathbb{D}},d\nu_q) \ar[d]^{\mathcal{L}^{\mathrm{radial}}_k}
\ar[r]^{\mathscr{I}} &
L^2(q^{-2\mathbb{Z}^n_+} ) \ar[d]^{\widetilde{\mathcal{L}_k}}\\
L^2(\Delta_{\mathbb{D}},d\nu_q) \ar[r]^{\mathscr{I}} &
L^2(q^{-2\mathbb{Z}^n_+}).} $$

\begin{lemma}\label{oper} The operators
$\mathcal{L}^{\mathrm{radial}}_1, \mathcal{L}^{\mathrm{radial}}_2, \ldots
\mathcal{L}^{\mathrm{radial}}_n$ in the Hilbert space
$L^2(\Delta_{\mathbb{D}},d\nu_q)$ are bounded self-adjoint and pairwise
commuting.
\end{lemma}
{\bf Proof.} The explicit formulas of Subsection \ref{Plansh_1} imply that the
operators $\Box_{u_i}$ are bounded for all $1 \leq i \leq n$ (unlike in the
classical case), so the same holds for $\widetilde{\mathcal{L}}_1,
\widetilde{\mathcal{L}}_2, \ldots \widetilde{\mathcal{L}}_n$. Moreover, it is easy
to see that the operators $\Box_{u_i}$ and $\Box_{u_j}$ commute for $1 \leq i <
j \leq n$ (as they are acting in different variables), so the operators
$\widetilde{\mathcal{L}_i},\widetilde{\mathcal{L}_j}$ commute for $1 \leq i < j
\leq n$ too.

Also, the operators $\Box_{u_i}$, $i=1,2,\ldots,n$, are symmetric, so they
are bounded self-adjoint operators in $L^2(q^{-2\mathbb{Z}^n_+})_q$. Thus,
$\widetilde{\mathcal{L}_i}$, $i=1,2,\ldots,n$, are pairwise commuting
bounded self-adjoint linear operators in $L^2(q^{-2\mathbb{Z}^n_+})_q$. As
the mapping $\mathscr{I}$ is isometric, the operators
$\mathcal{L}^{\mathrm{radial}}_1, \mathcal{L}^{\mathrm{radial}}_2, \ldots
\mathcal{L}^{\mathrm{radial}}_n$ are also bounded self-adjoint and pairwise
commuting.\hfill $\Box$

Using Proposition \ref{own_rk1}, one can easily show that the functions
$\Phi_{l_1}(u_1) \Phi_{l_2}(u_2) \ldots \Phi_{l_n}(u_n)$ on
$q^{-2\mathbb{Z}_+^n}$ are common generalized eigenfunctions of the
operators \eqref{raznvyr}. We will need the common eigenfunctions which are in
the image of the operator $\mathscr{I}$. It is easy to see that
$$\widetilde{\phi}_{l_1,l_2,\ldots,l_n}(u_1,u_2,\ldots,u_n) =
\sum\limits_{\sigma \in S_n} \mathrm{sign}(\sigma) \;
\Phi_{l_1}(u_{\sigma_1}) \Phi_{l_2}(u_{\sigma_2}) \ldots
\Phi_{l_n}(u_{\sigma_n}) \in \mathrm{Im} \mathscr{I}$$ are common
generalized eigenfunctions. 
Let
\begin{equation*}
\mathcal{R} = \{ (\rho_1,\rho_2,\ldots,\rho_n) \in [0,\pi/h]^n,
\quad \rho_1>\rho_2\>\ldots>\rho_n\}.
\end{equation*}

\begin{lemma}\label{oper_razl}
The pairwise commuting bounded self-adjoint operators
$\widetilde{\mathcal{L}_k}$, $k=1,2,\ldots,n$, are unitary equivalent to the
operators of multiplication by $$e_k( a(-\frac{1}{2}+i
\rho_1),a(-\frac{1}{2}+i \rho_2),\ldots,a(-\frac{1}{2}+i \rho_n)),\qquad
k=1,2,\ldots,n,$$ (respectively) in the Hilbert space $L^2(\mathcal{R},
(n!)\mathcal{N} (d\sigma)^n |_{\mathcal{R}})$. The unitary equivalence is
provided by the mapping
$$
\widetilde{\mathcal{U}}: \mathrm{Im} \;\mathscr{I} \to L^2(\mathcal{R}, (n!)
\mathcal{N} (d\sigma)^n |_{\mathcal{R}}),$$
$$
\widetilde{\mathcal{U}}:\; f(u_1,u_2,\ldots,u_n) \mapsto $$ $$
\widehat{f}(\rho_1,\rho_2,\ldots,\rho_n)= \mathcal{N} \int\limits_1^\infty \ldots
\int\limits_1^\infty \widetilde{\phi}_{-\frac{1}{2}+i \rho_1,-\frac{1}{2}+i \rho_2,
\ldots,-\frac{1}{2}+i \rho_n}(u) f(u) d_{q^{-2}} u_1 \ldots d_{q^{-2}} u_n.
$$
The inverse operator is
$$ \widetilde{\mathcal{U}}^{-1} : \widehat{f}(\rho_1,\rho_2,\ldots,\rho_n)
\mapsto $$ $$ \mathcal{N} \underbrace{\int\ldots\int}_\mathcal{R}
\widehat{f}(\rho_1,\rho_2, \ldots,\rho_n) \widetilde{\phi}_{-\frac{1}{2}+i
\rho_1, -\frac{1}{2}+i \rho_2,\ldots,-\frac{1}{2}+i \rho_n}(u) \; (n!)
d\sigma(\rho_1) \ldots d\sigma(\rho_n).$$
\end{lemma}
{\bf Proof.} This lemma follows from results of \itemii $ $ \ref{Plansh_1} and the
explicit formulas for the operators $\widetilde{\mathcal{L}_1}$,...,
$\widetilde{\mathcal{L}_n}$.
\begin{remark}
The last equalities define $\widetilde{\mathcal{U}}$ on a dense linear manifold of
the functions with finite support on the set $q^{-2\mathbb{Z}_+^n}$.
\end{remark}


Let us introduce the notation
\begin{equation}\label{phi_general_def}\Phi_{l_1,l_2,\ldots,l_n}(u)=
\frac{\sum\limits_{\sigma \in S_n} \mathrm{sign}(\sigma) \;
\Phi_{l_1}(u_{\sigma_1}) \Phi_{l_2}(u_{\sigma_2}) \ldots
\Phi_{l_n}(u_{\sigma_n})}{\Delta(u)}.\end{equation}
\begin{remark} (see Corollary \ref{det_formula} and Remark \ref{phi_capital}).
The spherical function $\varphi_\lambda(u)$, $\lambda \in \Lambda_n$ is equal
up to a multiplicative constant to $\Phi_{l_1,l_2,\ldots,l_n}(u)$, where
$l_i=(\lambda+\delta)_{i}\in \mathbb{Z}$.
\end{remark}

Using this lemma and the definition \eqref{isooper} of the operator
$\mathscr{I}$, one can easily obtain the following lemma.

\begin{lemma}\label{twoside_plansh_pre}
The pairwise commuting bounded self-adjoint operators
$\mathcal{L}^{\mathrm{radial}}_k$, $k=1,2,\ldots,n$ are unitary equivalent
to the operators of multiplication by $$e_k( a(-\frac{1}{2}+i
\rho_1),a(-\frac{1}{2}+i \rho_2),\ldots,a(-\frac{1}{2}+i \rho_n)),\qquad
k=1,2,\ldots,n$$ (respectively) in the Hilbert space $L^2(\mathcal{R},
(n!)\mathcal{N}(d\sigma)^n |_{\mathcal{R}})$. The unitary equivalence is
provided by the mapping
$$
\mathcal{U}: L^2(\Delta_{\mathbb{D}},d\nu_q) \to L^2(\mathcal{R},
(n!)\mathcal{N}(d\sigma)^n |_{\mathcal{R}}),$$
$$
\mathcal{U}:\; f(u) \mapsto \widehat{f}(\rho_1,\rho_2, \ldots,\rho_n)=
\int\limits_{\Delta_{\mathbb{D}}} \Phi_{-\frac{1}{2}+i
\rho_1,-\frac{1}{2}+i \rho_2,\ldots,-\frac{1}{2}+i \rho_n}(u) f(u)
d\nu_q(u),
$$
where the measure $d\nu_q(u)$ is defined in \eqref{measure_n}.

The inverse operator is
$$ \mathcal{U}^{-1} : \widehat{f}(\rho_1,\rho_2,\ldots,\rho_n)
\mapsto $$ $$ \int\limits_\mathcal{R} \widehat{f}(\rho_1,\rho_2,
\ldots,\rho_n) \Phi_{-\frac{1}{2}+i \rho_1, -\frac{1}{2}+i
\rho_2,\ldots,-\frac{1}{2}+i \rho_n}(u) \; (n!) \mathcal{N}\;
d\sigma(\rho_1) \ldots d\sigma(\rho_n).$$
\end{lemma}

\begin{remark}
The last equalities define $\mathcal{U}$ on a dense linear manifold of functions
with finite support on the set $\Delta_{\mathbb{D}}$.
\end{remark}

\begin{lemma}\label{f0_im}
\begin{equation}\label{fo_eq}\mathcal{U} f_0 = \mathcal{N} \Delta(q^{-2\delta})^{-1}
\left( \prod\limits_{j=0}^{n-1} \frac{(q^{-2j};q^2)_j}{(q^2;q^2)_j^2}
q^{(j+1)^2-1} \right) \prod\limits_{1 \leq k < j \leq n} (q^{-2i\rho_j} +
q^{2i\rho_j} - q^{-2i\rho_k} - q^{2i\rho_k}),\end{equation} where the
constant $\mathcal{N}$ is defined in \eqref{N_norm}.
\end{lemma}
{\bf Proof.}
\begin{multline}\label{f0_left}
(\mathcal{U} f_0 )(\rho_1, \rho_2, \ldots, \rho_n) = \mathcal{N}
\Phi_{-\frac{1}{2}+i \rho_1, -\frac{1}{2}+i \rho_2,\ldots,-\frac{1}{2}+i
\rho_n}(1, q^{-2}, \ldots, q^{-2(n-1)}) = \\
\mathcal{N} \Delta(q^{-2\delta})^{-1} \sum\limits_{\sigma \in S_n}
\mathrm{sign}(\sigma) \; \Phi_{-\frac{1}{2}+i \rho_1}(1)
\Phi_{-\frac{1}{2}+i \rho_2}(q^{-2}) \ldots \Phi_{-\frac{1}{2}+i
\rho_n}(q^{-2(n-1)}) .\end{multline}

It can be verified that the last expression is a polynomial in the variables
$q^{i \rho_1}+q^{-i \rho_1}, \ldots, {q^{i \rho_n}+q^{-i \rho_n}}$. It is
antisymmetric, so
\begin{equation}\label{fo_right}
\prod\limits_{1 \leq k < j \leq n} (q^{-2i\rho_j} + q^{2i\rho_j} -
q^{-2i\rho_k} - q^{2i\rho_k})
\end{equation}
is a factor of \eqref{f0_left}. One can compare the degrees of the
polynomials in the right-hand sides of \eqref{f0_left} and \eqref{fo_right}
as the elements of the graded algebra $\mathbb{C}[q^{i \rho_1}+q^{-i
\rho_1},q^{i \rho_2}+q^{-i \rho_2}, \ldots, q^{i \rho_n}+q^{-i \rho_n}]$.
The degree of the polynomial (\ref{fo_right}) is $\frac{n(n-1)}{2}$. Since
$$\Phi_{-\frac{1}{2}+i
\rho}(q^{-2k}) = {_{3}\Phi_{2}}\left(
\begin{array}{c}
q^{1+i\rho},q^{1-i\rho},q^{-2k}\\ q^2,0
\end{array}
;q^2,q^2\right) = $$ $$\sum\limits_{j=0}^k \frac{(q^{1+i\rho};q^2)_j\,
(q^{1-i\rho};q^2)_j \,(q^{-2k};q^2)_j \,q^{2j}}{(q^2;q^2)_j^2} ,$$ then the
degree of $\mathcal{U} f_0$ is $\frac{n(n-1)}{2}$, and it proves \eqref{fo_eq} up
to a constant. This constant can be found by comparing the highest monomial
coefficients in the lexicographic order.
 \hfill $\Box$

Denote
\begin{multline}\label{kappa}
\kappa(\rho_1,\rho_2,\ldots,\rho_n) = \\ \mathcal{N}
\Delta(q^{-2\delta})^{-1} \left( \prod\limits_{j=0}^{n-1}
\frac{(q^{-2j};q^2)_j}{(q^2;q^2)_j^2} q^{(j+1)^2-1} \right) \prod\limits_{1
\leq k < j \leq n} (q^{-2i\rho_j} + q^{2i\rho_j} - q^{-2i\rho_k} -
q^{2i\rho_k}). \end{multline}

Notice that the function $\kappa(\rho_1,\rho_2,\ldots,\rho_n)$ is positive on
$\mathcal{R}$. Consider the operator
\begin{equation*}\label{U_F}
\mathcal{F} = \frac{1}{\kappa(\rho_1,\rho_2,\ldots,\rho_n)}\; \mathcal{U}
\end{equation*}
and the measure
\begin{equation}\label{Sigma_sigma}
d\Sigma(\rho_1,\rho_2,\ldots,\rho_n) = \kappa(\rho_1,\rho_2,\ldots,\rho_n)^2
(n!)\mathcal{N} (d\sigma(\rho_1) \ldots d\sigma(\rho_n))|_{\mathcal{R}}
\end{equation}
on the set $\mathcal{R}$ (the constant $\mathcal{N}$ is defined in
\eqref{N_norm}).

The following proposition is the consequence of Lemmas
\ref{twoside_plansh_pre} and \ref{f0_im}.
\begin{proposition}\label{twoside_plansh}
The pairwise commuting bounded self-adjoint operators
$\mathcal{L}^{\mathrm{radial}}_k$, $k=1,2,\ldots,n$, are unitary equivalent to
the operators of multiplication by
$$\frac{e_k( a(-\frac{1}{2}+i \rho_1),a(-\frac{1}{2}+i
\rho_2),\ldots,a(-\frac{1}{2}+i
\rho_n))}{\kappa(\rho_1,\rho_2,\ldots,\rho_n)},\qquad k=1,2,\ldots,n$$
(respectively) in the Hilbert space $L^2(\mathcal{R}, d\Sigma)$. The unitary
equivalence is provided by the mapping
$$
\mathcal{F}: L^2(\Delta_{\mathbb{D}},d\nu_q) \to L^2(\mathcal{R},
d\Sigma),$$
\begin{multline}
\mathcal{F}:\; f(u) \mapsto \widehat{f}(\rho_1,\rho_2, \ldots,\rho_n)= \\
\frac{1}{\kappa(\rho_1,\rho_2,\ldots,\rho_n)}
\int\limits_{\Delta_{\mathbb{D}}} \Phi_{-\frac{1}{2}+i
\rho_1,-\frac{1}{2}+i \rho_2,\ldots,-\frac{1}{2}+i \rho_n}(u) f(u)
d\nu_q(u),
\end{multline}
where $\Phi_{l_1,l_2,\ldots,l_n}(u)$ are defined in \eqref{phi_general_def},
and the measure $d\nu_q(u)$ is defined in \eqref{measure_n}.

The inverse mapping is
$$ \mathcal{F}^{-1} : \widehat{f}(\rho_1,\rho_2,\ldots,\rho_n)
\mapsto $$ $$ \int\limits_\mathcal{R} \widehat{f}(\rho_1,\rho_2,
\ldots,\rho_n) \Phi_{-\frac{1}{2}+i \rho_1, -\frac{1}{2}+i
\rho_2,\ldots,-\frac{1}{2}+i \rho_n}(u) \;
d\Sigma(\rho_1,\rho_2,\ldots,\rho_n).$$
\end{proposition}

 \begin{remark}\label{fo_im} The cyclic vector $f_0
\in L^2(\Delta_{\mathbb{D}},d\nu_q)$ is mapped into $1 \in L^2(\mathcal{R},
d\Sigma)$ by $\mathcal{F}$.
\end{remark}

\begin{remark}\label{independ_var}
For the convenience we use the variables $\rho_1,\rho_2,\ldots,\rho_n$ in
the image of $\mathcal{F}$. Notice that if we change the variables $$z_k =
\frac{e_k( a(-\frac{1}{2}+i \rho_1),a(-\frac{1}{2}+i
\rho_2),\ldots,a(-\frac{1}{2}+i
\rho_n))}{\kappa(\rho_1,\rho_2,\ldots,\rho_n)}$$ we get that $\mathcal{F}$
maps $\mathcal{L}^{\mathrm{radial}}_k$ into the operator of multiplication
by the independent variable $z_k$.
\end{remark}

The measure $d\Sigma$ on $\mathcal{R}$ defined in \eqref{Sigma_sigma} is a
sought-for radial part of the Plancherel measure.

\end{document}